\documentclass[preprint,10pt]{elsarticle}

\usepackage{amssymb,graphicx,epsfig,amsmath,amsfonts}
\usepackage{color,algorithm,algpseudocode}

\newcommand{\Rr}{{\mathbb{R}}}
\newcommand{\Nn}{{\mathbb{N}}}

\def\a2{{\alpha/2}}

\def\ga{{\gamma}}

\def\bfq{{\mathbf{q}}}

\def\bfg{{\mathbf{g}}}

\def\bfz{{\mathbf{z}}}

\def\a2{{\alpha/2}}

\def\GG{{\Gamma}}

\def\bfzeta{\mbox{\boldmath$\zeta$}}
\def\qed{\hfill$\Box$}
\def\cR{{\cal{R}}}
\def\cP{{\cal{P}}}
\def\cU{{\cal{U}}}

\def\llangle{\left\langle}
\def\rrangle{\right\rangle}

\newtheorem{theorem}{Theorem}[section]
\newtheorem{remark}{Remark}[section]
\newtheorem{proposition}{Proposition}[section]
\newtheorem{lemma}{Lemma}[section]

\begin{document}

\begin{frontmatter}

\title{A corrected spectral method for Sturm-Liouville problems with unbounded potential at one endpoint.}

\author[auto]{Cecilia Magherini}
\ead{cecilia.magherini@unipi.it}
\address[auto]{Dipartimento di Matematica, Universit\`{a} di Pisa, I-56127 Pisa (Italy).}
\begin{abstract}
In this paper, we shall derive a spectral matrix method for the approximation of the eigenvalues of (weakly) regular and singular
Sturm-Liouville problems in normal form with an unbounded potential at the left endpoint. The method is obtained 
by using a Galerkin approach with an approximation of the eigenfunctions given by suitable combinations of Legendre polynomials.
We will study the errors in the eigenvalue estimates for problems with unsmooth eigenfunctions  in 
proximity of the left endpoint. The results of this analysis will be then used conveniently to determine low-cost and 
effective procedures for the computation of corrected numerical eigenvalues. Finally,  we shall present and discuss 
the results of several numerical experiments which  confirm the effectiveness of the approach.
\end{abstract}

\end{frontmatter}

\section{Introduction}
The direct Sturm-Liouville problem (SLP) in normal form with separated boundary conditions is given by
\begin{eqnarray}\label{slp}
-y''(x) + q(x) y(x) &=& \lambda y(x), \qquad  x \in (a,b)\>, \\
\alpha_a y(a) + \beta_a y'(a) &=& 0\>, \qquad \quad\alpha_a^2+\beta_a^2 \neq 0,\label{bca}\\
\alpha_b y(b) + \beta_b y'(b) &=& 0\>, \qquad \quad\alpha_b^2+\beta_b^2 \neq 0,\label{bcb}
\end{eqnarray}
where the potential $q,$ the domain $(a,b)$ and the 
coefficients 
$\alpha_{a},\beta_a,\alpha_b,\beta_b$ 
represent the data of the problem while the 
unknowns are the eigenvalues $\lambda$ and the corresponding eigenfunctions $y.$  
It is surely a classical problem that has been extensively studied both from the theoretical and from the 
numerical point of views. Many numerical schemes are nowadays available for its solution which can be 
subdivided into two main families: matrix methods and shooting techniques \cite{Pryce}.
A number of well-established numerical codes that are able to solve regular problems as well as many singular 
ones had been developed and are freely available for the scientific community. Among them we surely mention 
the MATSLISE2 \cite{Led}, the SLEDGE \cite{Sledge} and the SLEIGN2 \cite{Sleign} codes, 
 but there are many others.
In spite of this, we think that the approach of deriving matrix schemes by using spectral methods, instead of 
finite difference or element ones, deserves further insights and this is the topic of the present paper.
In particular, we will consider the case of a bounded domain which, without loss of generality, we assume to be
\begin{equation} \label{ab}
(a,b) \equiv (-1,1), 
\end{equation}
and study problems with 
\begin{equation} \label{vx}
 q(x) = f(x) + \frac{g(x)}{(1+x)^\ga},
\end{equation}
where $f$ and $g$ are analytic functions inside and on a Bernstein ellipse containing $[-1,1],$
and $\ga\ge 0.$ Problems of this type have many applications in physics (cf., for
example, \cite{Cal64, Cas50, UHH81, VW54}). We recall that if $q\in L_1(-1,1)$ 
then the problem is regular (sometimes called weakly regular if $\ga \in (0,1)$ and $g(-1) \neq 0$) 
otherwise it is singular.
More precisely, the singular left endpoint is of type \cite{AF85,DS,Kod,Tit}
\begin{enumerate}
 \item Limit-Circle (LC) and nonoscillatory  if $\ga\in [1,2)$ and $g(-1) \neq 0$ or $\ga=2$ and  $g(-1)\in [-1/4,3/4)\setminus\{0\};$
 \item Limit-Circle and oscillatory  if $\ga=2$ and $g(-1)<-1/4$ or $\ga >2$ and $g(-1) <0;$
 \item Limit-Point (LP)  and nonoscillatory  if $\ga=2$ and $g(-1)\ge 3/4$ or $\ga >2$ and $g(-1) >0.$
\end{enumerate} 
In the sequel, we will always exclude the case of an oscillatory endpoint which is definitely 
much more difficult to be treated numerically and we will assume $\ga \in [0,2]$
leaving the generalization to larger values of $\ga$ to future investigation.
Concerning the boundary condition at $x=-1,$
which is required in the LC case, we will consider 
the Friedrichs one, namely in this context the Dirichlet condition. 
This means that for any $\lambda$ we select the principal solution of the 
equation which is frequently the most significant in the applications, \cite{NZ,Pryce}.
We recall that the Dirichlet condition is the only possible one in the LP case.\\
Under these assumptions, it is known that the eigenvalues of (\ref{slp})--(\ref{vx}) are real, simple and that they can be ordered 
as an increasing sequence tending to infinity. We will number them starting from index $k=1,$ i.e. we will call
$$\left\{\lambda_1 < \lambda_2 < \lambda_3< \ldots\right\}$$
the exact spectrum of (\ref{slp})--(\ref{vx}).\\
Talking about numerical methods based on shooting techniques, the approach used most frequently for solving 
this type of problems is based on the selection of a layer, namely
(\ref{slp}) is usually solved over $(-1+\epsilon,1).$ In particular, suitable algorithms 
have been studied for an adaptive selection of $\epsilon$ and for the computation of the  
condition to be imposed at $x=-1+\epsilon$ (see, for example, \cite{Sleign,LedTh} for further details).\\
Concerning classical matrix methods, it is clear that they can be employed if based on a discretization of the differential equation and 
of the boundary conditions which do not require the evaluation of $q,$ or even of its derivative,  at the left endpoint. 
This is the case, for example, of the classical three-point formula. The main difficulty with 
this approach may be that of an order reduction caused by the fact that the derivative, of suitable order, of the eigenfunctions    
may be unbounded as $x$ approaches $-1.$\\

In this paper, we shall derive a matrix method by using a spectral Galerkin approach 
based on Legendre polynomials.  Before proceeding, it must be said that if the problem is subject to 
Dirichlet boundary conditions at both endpoints then  
schemes based on spectral (collocation) methods 
which use orthogonal polynomials or sinc functions are already available in the literature (see, for example, \cite{CS,El,GSS,JLB}). 
These are surely effective methods to be employed whenever the eigenfunctions are sufficiently smooth.  
Our main purpose is therefore that of treating general boundary conditions and to get accurate approximations of the 
eigenvalues even in the case where the eigenfunctions are not so regular.\\

The remaining part of this paper is organized as follows. In Section~\ref{sec2}, we describe the approach, derive the generalized 
eigenvalue problem which discretize the continuous one and discuss how the entries of the matrices involved can be 
computed efficiently. Section~\ref{sec3} is devoted to the analysis of the errors in the resulting eigenvalue
approximations for problems with unbounded potential at the left endpoint. Moreover, in the same section we 
shall derive low cost and effective procedures for the computation of corrected numerical eigenvalues. Finally, the results of several
numerical experiments are reported and discussed in Section~\ref{sec4}. 

\section{Spectral Legendre-Galerkin method}\label{sec2}
Let $\Pi_{N+1}$ be the space of polynomials of maximum degree $N+1,$ 
for a fixed $N\in \Nn,$ and let
\begin{eqnarray}
{\cal S}_N &\equiv& \left\{  r \in \Pi_{N+1}: \quad 
 \alpha_{a} \,r(-1) + \beta_a\, r(-1) = \alpha_{b}\, r(1)+ \beta_b \,r(1) =0 \right\}\label{SN}\\
       &\equiv& \mbox{span} \left(\cR_0,\cR_1,\ldots,\cR_{N-1}\right). \label{spanSN}
             \end{eqnarray}
We look for an approximation of an eigenfunction $y$ of the following type
\begin{equation}\label{zN}
z_N(x) = \sum_{n=0}^{N-1} \zeta_{n,N} \cR_n(x) \approx y(x)
\end{equation}
where the coefficients $\zeta_{n,N}$ and the numerical eigenvalue $\lambda^{(N)}$ are determined by imposing 
\begin{equation}\label{wf}
\sum_{n=0}^{N-1} \left\langle \cR_{m},-\cR_n'' + (q-\lambda^{(N)}) \cR_n\right\rangle \zeta_{n,N}=0, \quad \mbox{for each } m=0,\ldots,N-1.
\end{equation}
Here $\langle \cdot,\cdot \rangle$ is the standard inner product in $L_2([-1,1]),$ i.e.
$$ \langle u, v\rangle = \int_{-1}^1 u(x)v(x) dx, \qquad u,v \in L_2([-1,1]),$$
which is naturally suggested by the Liouville normal form of the SLP we are studying.
We can write (\ref{wf})  as the following generalized eigenvalue problem
\begin{equation}\label{geneig}
\left(A_N + Q_N\right)\bfzeta_N = \lambda^{(N)} B_N \bfzeta_N  
\end{equation}
where $\bfzeta_N = \left(\zeta_{0N}, \ldots,\zeta_{N-1,N}\right)^T,$ 
\begin{equation}\label{ABQ}
A_N = \left(a_{mn}\right), \quad B_N = \left(b_{mn}\right), \quad 
Q_N = \left(q_{mn}\right), \quad m,n=0,\ldots,N-1,
\end{equation}
with
\begin{eqnarray}
\label{abmn}
a_{mn} &=&  -\langle \cR_m, \cR_n''\rangle, \quad b_{mn} = \langle \cR_m, \cR_n\rangle, \\
q_{mn} &=& \langle \cR_m, f\,\cR_n\rangle + \langle \cR_m,(1+x)^{-\ga} g \cR_n\rangle \equiv f_{mn} + g_{mn}. \label{qzmn}
\end{eqnarray}
The matrices $B_N$ and $Q_N$ are clearly symmetric. The same property holds for $A_N$ thanks to
the well-known Green's identity 
\begin{equation}\label{gid}
\langle v,u''\rangle -\langle u,v''\rangle   = \left[u'(x)v(x) -u(x)v'(x) \right]_{-1}^1,
\end{equation}
by which one gets, for each  $\cR_n$ and $\cR_m \in {\cal{S}}_N,$
$$\langle \cR_m,\cR_n''\rangle = \langle \cR_n,\cR_m''\rangle + \left[\cR_n'(x)\cR_m(x) -\cR_n(x)\cR_m'(x) \right]_{-1}^1 =
\langle \cR_n,\cR_m''\rangle. $$

Clearly the effectiveness of the procedure is strictly connected to the choice of the basis functions.
The main criterion we have considered is the computational cost of the method which is essentially 
determined by the calculation of the coefficient matrices and by the solution of (\ref{geneig}).
This suggests to use suitable combinations of the classical Legendre polynomials as described  in the 
next subsection.

\subsection{Basis functions}
As done in \cite{Shen94}, we look for a basis function $\cR_{n}$ of the following form 
\begin{equation}\label{Rnp2}
 \cR_{n}(x) = \xi_{n} \cP_n(x) +  \eta_{n} \cP_{n+1}(x) +  \theta_{n} \cP_{n+2}(x) 
\end{equation}
where $\cP_j$ is the Legendre polynomial of degree $j$ for which
it is known that \cite{Abram}
\begin{equation}\label{valP}
\cP_j(1) = (-1)^j \cP_j(-1) = 1, \qquad \cP'_j(1) = (-1)^{j-1}\cP'_j(-1) = j(j+1)/2.
\end{equation}
Therefore, with some computations one gets that $\cR_n \in {\cal S}_N,$ see (\ref{SN})-(\ref{spanSN}), 
if and only if $\left(\xi_n,\eta_n,\theta_n\right)^T$ belongs to the kernel of 
$$
V_{n} = \left(\begin{array}{r|r|r}
     \alpha_a - \frac{n(n+1)}2\beta_a & -\alpha_a + \frac{(n+1)(n+2)}2\beta_a & \alpha_a - \frac{(n+2)(n+3)}2\beta_a\\
     &&\vspace{-2mm}\\
     \hline
     &&\vspace{-2mm}\\
     \alpha_b + \frac{n(n+1)}2\beta_b &  \alpha_b + \frac{(n+1)(n+2)}2\beta_b & \alpha_b + \frac{(n+2)(n+3)}2\beta_b
\end{array}\right).
$$

We must now distinguish the following two possibilities:
\begin{enumerate}
 \item $\alpha_a\beta_b + \alpha_b \beta_a = 0,$ i.e. problems subject to symmetric BCs. In this case it is natural to
   set $\eta_{n} = 0$ so that $\cR_{n}$ is an even or an odd function, depending on $n,$ which
   implies that if (\ref{bca}) holds true then (\ref{bcb}) is verified automatically. In this way, one obtains
   the following system of equations
   \begin{eqnarray*}
    \left\{\begin{array}{r}
      \left(\alpha_a - \frac{n(n+1)}2\beta_a\right) \xi_{n} +  \left(\alpha_a - \frac{(n+2)(n+3)}2\beta_a\right) \theta_{n} = 0,\\
      \eta_n = 0,
     \end{array}
   \right.
   \end{eqnarray*}
   whose general solution can be written as
   \begin{eqnarray} \label{xins}
   \xi_n &=& \,-\nu_n \left(\alpha_a - \frac{(n+2)(n+3)}2\beta_a\right),\\
   \eta_n &=& 0,\label{etans}\\
   \theta_n &=& \,\,\nu_n \left(\alpha_a - \frac{n(n+1)}2\beta_a\right), \label{tetans}
   \end{eqnarray}
   where $\nu_n\neq 0$ is a free parameter;
   
 \item $\alpha_a\beta_b + \alpha_b \beta_a \neq  0.$ From the previous considerations, one deduces that $\eta_{n}$ must be different from zero.
  Using the Matlab notation, this is confirmed by the fact that 
  $$\mbox{det}\left(V_n(\,:\,,[1~~3])\right)=(2n+3)(\alpha_a\beta_b+\alpha_b\beta_a)\neq 0, \qquad \mbox{for each } n.$$
 Hence, if we let as before $\nu_n\neq 0$ be a free parameter then we got
\begin{small}
\begin{eqnarray} 
\xi_n &=& -\nu_n \,\mbox{det}\left(V_n(\,:\,,[2~~3])\right)\nonumber \\
&=& \nu_n \,\left[2\alpha_a\alpha_b +(n+2)^2\left(\alpha_a\beta_b-\alpha_b\beta_a-\frac{(n+1)(n+3)}2\beta_a\beta_b\right)\right],\label{xinns}\\
\eta_n &=& \nu_n \,\mbox{det}\left(V_n(:,[1~~3])\right) 
= \nu_n\,(2n+3)(\alpha_a\beta_b+\alpha_b\beta_a), \label{etanns}\\
\theta_n &=& -\nu_n \,\mbox{det}\left(V_n(:,[1~~2])\right)\nonumber \\
&=&  -\nu_n \,\left[2\alpha_a\alpha_b +(n+1)^2\left(\alpha_a\beta_b-\alpha_b\beta_a-\frac{n(n+2)}2\beta_a\beta_b\right)\right].\label{tetanns}
\end{eqnarray}
\end{small}
\end{enumerate}

Concerning the parameters $\left\{\nu_n\right\}_{n\in\Nn_0}$  we decided 
to establish a criterion for their selection in order to obtain basis functions 
independent of the scaling of $(\alpha_a,\beta_a)$ and/or $(\alpha_b,\beta_b).$
Now, a natural choice would have been that of choosing $\nu_n$  so that $\|\cR_n\|_2=1$ for each $n$ since we are working in $L_2.$
Nevertheless, we preferred not to proceed in this way to avoid the computation of square roots at least at this level.
As an alternative to $\|\cR_n\|_2$ we considered its uniform norm which is not known in closed form but it satisfies
$\|\cR_n\|_\infty \leq 3 \|\left(\xi_n,\eta_n,\theta_n\right)^T\|_\infty$.
We thus applied the following criterion
$$
\mbox{for each $n$ let $\nu_n$ be such that }  \|\left(\xi_n,\eta_n,\theta_n\right)^T\|_\infty = 1  \mbox{ and } \xi_n \ge 0.
$$
The so-obtained coefficients for the four problems subject to natural BCs and for two general ones are listed in Table~\ref{tabR}.
The unspecified values for the last BCs, of Dirichlet-Robin type, are 
$\xi_0=2/3,\,\eta_0=1$ and $\theta_0=1/3.$\\

\begin{table}
 \caption{Coefficients $\xi_{n}, \eta_{n}$ and $\theta_{n}$ for some BCs.}\label{tabR}
\begin{center}
\begin{tabular}{|c|c||c|c|c|}
\hline
&&&&\vspace{-2mm}\\
BCs & Name & $\xi_{n}$ & $\eta_{n}$ & $\theta_{n}$\\
&&&&\vspace{-2mm}\\
\hline
\hline
&&&&\vspace{-3mm}\\
$y(\pm 1)=0$
& \small{Dirichlet--Dirichlet}& 1 & 0 & $-1$\\
&&&&\vspace{-3mm}\\
\hline
&&&&\vspace{-3mm}\\
$y'(\pm 1)=0$
& \small{Neumann--Neumann} & 1 & 0 & $-\frac{n(n+1)}{(n+2)(n+3)}$\\
&&&&\vspace{-3mm}\\
\hline
&&&&\vspace{-3mm}\\
$\begin{array}{l} y(-1)=0 \\ \,\,\,y'(1)=0\end{array}$
& \small{Dirichlet-Neumann} &  $1$ & $\frac{(2n+3)}{(n+2)^2}$ &  $-\left(\frac{n+1}{n+2}\right)^2$\\
&&&&\vspace{-3mm}\\
\hline
&&&&\vspace{-3mm}\\
$\begin{array}{l} y'(-1)=0 \\ \quad y(1)=0\end{array}$
& \small{Neumann-Dirichlet} &  $1$ & $-\frac{(2n+3)}{(n+2)^2}$ &  $-\left(\frac{n+1}{n+2}\right)^2$\\
&&&&\vspace{-3mm}\\
\hline
&&&&\vspace{-3mm}\\
$\begin{array}{l}
  y(-1)\,=\,-y'(-1)\\
  \,\,\,\,y(1)\,=\,\,\,\,y'(1)\\
\end{array}$
& \small{$\begin{array}{c} \mbox{Robin-Robin}\\ \mbox{symmetric}\end{array}$} &  $1$ & $0$ &  $-\,\frac{n(n+1)-2}{(n+2)(n+3)-2}$\\
&&&&\vspace{-3mm}\\
\hline
&&&&\vspace{-3mm}\\
$\begin{array}{l}
 y(-1) = 0\\
 \,\,\,y(1)\,=\,y'(1)
\end{array}$
& \small{$\begin{array}{c} \mbox{Dirichlet-Robin}\\  n \ge 1\end{array}$} & $1$ & $\frac{2n+3}{(n+2)^2-2}$ & $-\, \frac{(n+1)^2-2}{(n+2)^2-2}$ \\
\hline
\end{tabular}
\end{center}
\end{table}
\noindent Finally, for later reference, it is important to underline the fact that, as soon as $n$ is sufficiently large, 
we always got
\begin{eqnarray}
\xi_n &=& 1,\label{xias}\\
\theta_n &=& -1 + O(n^{-1}),\label{thas}\\
\label{etas}
\eta_n &=& \left\{\begin{array}{ll}
                0, & \mbox{ if $\alpha_a\beta_b+\alpha_b\beta_a = 0,$}\\
                O(n^{-1}) & \mbox{ if $\alpha_a\beta_b+\alpha_b\beta_a \neq 0$ and $\beta_a\beta_b =0,$}\\
                O(n^{-3}) & \mbox{ if $\alpha_a\beta_b+\alpha_b\beta_a \neq 0$ and $\beta_a\beta_b \neq 0.$}
                \end{array}\right.          
\end{eqnarray}

\subsection{The matrices $A_N$ and $B_N$}
In this section we are going to show that the entries of $A_N$ and $B_N$ in (\ref{ABQ})-(\ref{abmn})
can be determined analytically thanks predominantly to the orthogonality of the Legendre polynomials with respect to
the standard inner product. \\

Concerning the first matrix, one immediately gets that $a_{mn}=0$ for each $m>n$
since $\cR_m$ is orthogonal to any polynomial in $\Pi_{m-1},$ see (\ref{Rnp2}). Consequently,  
$A_N=A_N^T$ is diagonal with diagonal entries 
\begin{eqnarray}\nonumber
a_{nn} &=& -\xi_{n} \langle \cP_n,\cR_n''\rangle - \eta_n \langle \cP_{n+1},\cR_n''\rangle - \theta_n \langle \cP_{n+2},\cR_{n}''\rangle\\
    &=& -\xi_n \langle \cP_n,\cR_n''\rangle = -\xi_n\theta_n \langle \cP_n,\cP_{n+2}''\rangle \nonumber\\
    &=& -\xi_n \theta_n\left[\cP_n(x)\cP_{n+2}'(x) - \cP_n'(x)\cP_{n+2}(x) \right]_{-1}^1  \nonumber\\
    &=& - 2 (2n+3) \xi_n\theta_n\label{amm1}
\end{eqnarray}
where for the last two equalities we used (\ref{gid}) and (\ref{valP}). For later convenience, we remark that 
independently of the BCs, $a_{nn}$ satisfies 
\begin{equation}\label{annas}
a_{nn} = 4\left(n+\frac{3}{2}\right) + O\left(n^{-1}\right), \quad n \gg 1.  
\end{equation}

Regarding $B_N,$ it is not too difficult to verify that it is pentadiagonal. In more detail, if we let 
$$ \hat{b}_{n}  = \langle \cP_n,\cP_n\rangle = 2/(2n+1),$$
\begin{equation} \label{R}
 \hat{B}_N = \left(\begin{array}{ccc} \hat{b}_{0}\\ &\ddots \\ &&\hat{b}_{N+1}\end{array}\right), \quad
  R_N = \left(\begin{array}{ccccc}
       \xi_0\\
       \eta_0 & \ddots\\
       \theta_0 & \ddots & \ddots\\
       & \ddots & \ddots & \xi_{N-1}\\
       && \ddots & \eta_{N-1}\\
       &&& \theta_{N-1}
\end{array}\right),
\end{equation}
then we get
\begin{equation}\label{Bfact}
B_N = R_N^T \, \hat{B}_N R_N.
\end{equation}

\subsection{The matrix $Q_N.$}\label{secQ}
Let us consider first of all the case of a regular problem with $\gamma <1.$ 
From (\ref{qzmn}), one obtains that $Q_N$
admits a factorization similar to the one just given for $B_N.$ Specifically
\begin{eqnarray} \label{Qfact}
Q_N &=& R_N^T \, \hat{Q}_N R_N, \\
\label{hFG}
\hat{Q}_N &=& \left(\hat{q}_{mn}\right) = \left(\hat{f}_{mn}+ \hat{g}_{mn}\right) \equiv \hat{F}_N + \hat{G}_N \in \Rr^{(N+2)\times(N+2)},
\end{eqnarray}
where $R_N$ is defined in (\ref{R}) and 
\begin{eqnarray}
\hat{f}_{mn} &=&  \int_{-1}^{1} f(x) \cP_m(x) \cP_n(x) dx, \label{hfmn}\\
\hat{g}_{mn} &=&  \int_{-1}^{1} (1+x)^{-\gamma} g(x) \cP_m(x) \cP_n(x) dx. \label{hgmn}
\end{eqnarray}
We recall that the Legendre polynomials obey the recurrence relation
\begin{eqnarray*}\nonumber
\cP_{-1} (x) &\equiv& 0, \qquad \qquad 
\cP_0 (x) \equiv 1,\\
\cP_{n+1}(x) &=& \frac{2n+1}{n+1} x\, \cP_n(x) - \frac{n}{n+1} \cP_{n-1}(x), \quad n \ge 0. 
\end{eqnarray*}
This allows to prove the following result.

\begin{proposition}\label{propGFinf} Let $q\in L_1([-1,1])$ and, see (\ref{hFG}), let
\begin{eqnarray*}
\hat{\bfq}_n &\equiv& \left(\hat{q}_{0n}, \hat{q}_{1n}, \ldots \right )^T  \in \ell_\infty, \quad n \ge -1,
\end{eqnarray*}
with $\hat{\bfq}_{-1}$ the zero sequence.
If we define the linear tridiagonal operator $ \bfz \in \ell_\infty \mapsto \mathcal{H}\,\bfz \in \ell_\infty$ where
\begin{eqnarray}\nonumber
\mathcal{H} &=&
\left( \begin{array}{cccccc}
            0 & h_{01}\\
            h_{10} & 0 & h_{12}\\
            &h_{21} & 0 & h_{23}\\
            && \ddots & \ddots & \ddots\\
            \end{array}\right), \quad \begin{array}{l} h_{m,m-1}=m/(2m+1), \\ \\ h_{m,m+1}=(m+1)/(2m+1),\end{array}    \label{H}
\end{eqnarray}
then we get 
\begin{equation} \label{recvn}
\hat{\bfq}_{n+1} = \frac{2n+1}{n+1}\, \mathcal{H}\hat{\bfq}_n -\frac{n}{n+1} \,\hat{\bfq}_{n-1},\qquad  n\ge 0.
\end{equation}
\end{proposition}
\underline{Proof}: see \cite[Propositions 1,2]{GM} with $\alpha=0.$\qed\\

The immediate consequence of this proposition is that the recurrence in 
(\ref{recvn}) permits to 
determine the entire matrix $\hat{Q}_N$ once $\hat{q}_{m0}=\hat{f}_{m0}+\hat{g}_{m0}$ 
have been computed for each $m=0,1,\ldots,2N+2.$ 
In fact, these values are sufficient to determine $\hat{q}_{m1}$ for each $m=0,1,\ldots,2N+1$
by using (\ref{recvn}) with $n=0$ and the fact that $\mathcal{H}$ is tridiagonal. At this point,
the application of (\ref{recvn}) with $n=1$ allows to compute $\hat{q}_{m2}$ for each $m=0,1,\ldots,2N,$ and so forth.
Clearly,  in the actual implementation, the symmetry of $\hat{Q}_N$ is taken into account.\\
In particular, the coefficients $\hat{f}_{m0}$ in (\ref{hfmn}) decay exponentially as $m$ increases due to the 
assumption that $f$ is analytical inside and over a Bernstein ellipses containing $[-1,1].$ In addition, we can use the
routine {\tt legcoeffs} of the well-established open-source software package {\tt Chebfun} \cite{Cf} 
to determine  the numerically significant values. Thus if $L+1$ is the length of the vector 
of Legendre coefficients of $f$ provided by such routine, we approximate upto machine precision
the symmetric matrix $\hat{F}_N$  with its banded portion with bandwidth $2L+1.$
\begin{remark} 
 It is important to stress that if $q$ is analytical, i.e. if $g\equiv 0,$ then  
 the generalized eigenvalue problem (\ref{geneig})  which discretize the SLP involves only sparse matrices. In particular,
 $B_N$ is symmetric positive definite and pentadiagonal while $A_N+Q_N = A_N+R_N^T \hat{F}_N R_N$ is symmetric with bandwith $2L+5.$
\end{remark}
Concerning the computation of the required entries of $\hat{\bfg}_{0}=\left(\hat{g}_{10},\hat{g}_{20},\ldots\right)^T,$  
see (\ref{hgmn}),  we applied arguments similar to the ones used in \cite[Propositions 2,3]{GM}. In detail,
recalling that by assumption $g$ is analytical inside and over a Bernstein ellipse 
containing $[-1,1]$ too, the operator
$ g(\mathcal{H}) = \sum_{\ell=0}^{+\infty} \frac{\langle g,\cP_\ell \rangle}{\langle \cP_\ell,\cP_\ell \rangle}\, \cP_\ell\left(\mathcal{H}\right)$
is well defined. Consequently, \cite[Proposition 2]{GM} and \cite[16.4 formula (2)]{EMOT} 
allow to get that if $\gamma<1$ then 
\begin{equation}\label{hg0}
\hat{\bfg}_0 = g(\mathcal{H}) \left(\begin{array}{c}\hat{g}_0\\ \hat{g}_1\\\vdots \end{array}\right), \quad
\hat{g}_m = \int_{-1}^{1} \frac{\cP_m(x)}{(1+x)^{\gamma}} dx =  \frac{(-1)^m\,2^{1-\gamma}\,(\gamma)_m}{(1-\gamma)_{m+1}}, 
\end{equation}
where $(t)_\ell$ is the Pochhammer symbol. We observe that $\hat{g}_m$ coincides with $\hat{g}_{m0}$ if $g(x) \equiv 1$.
In the actual implementation, we proceed as follows: we get a polynomial approximation of $g$ by transforming it in a 
{\tt Chebfun} function, which is accurate up to machine precision, then we apply 
the previous formula to compute the first $2N+2$ entries of $\hat{\bfg}_0.$\\

Let us now discuss the case of singular problems, namely how we determine $Q_N$ if $\gamma \in [1,2].$
We recall that the corresponding Friedrichs  boundary condition at the singular endpoint is $y(-1)=0.$
The basis functions have therefore a root at  $x=-1$ and we need to highlight this fact.
This is done in the following proposition.
\begin{proposition}\label{propU} If $\beta_a=0$ and
$\cP_\ell^{(0,1)}$ is the Jacobi polynomial of degree $\ell$ with weighting function $\omega(x) = (1+x),$
then $\cR_n$ in (\ref{Rnp2}) can be written as 
\begin{equation}\label{Rndm1}
 \cR_n(x) = (1+x) \cU_n(x) = (1+x) \left(\xi_n \cP_n^{(0,1)}(x) + \theta_n \cP_{n+1}^{(0,1)}(x)\right).
\end{equation}
\end{proposition}
\underline{Proof:} The first equality is evident with $\cU_n\in \Pi_{n+1}$ since $\cR_n \in \Pi_{n+2}.$
In addition  
$$ \int_{-1}^1 (1+x) \cU_{n}(x)\,v(x) dx = <\cR_n,v> =  0, \quad \mbox{ for each } v\in \Pi_{n-1}.$$
This implies that there exist suitable $\tilde{\xi}_n$ and $\tilde{\theta}_n$ such that  
$$\cU_{n}(x) = \tilde{\xi}_n \cP_n^{(0,1)}(x) + \tilde{\theta}_n \cP_{n+1}^{(0,1)}(x).$$
It remains to prove that $\tilde{\theta}_n = \theta_n$ and $\tilde{\xi}_n=\xi_n.$
The former equality is an immediate consequence of the fact that 
$\cP_{n+2}$ and $\cP_{n+1}^{(0,1)}$ have the same leading coefficient, \cite{Abram}. 
Concerning the latter equality, it is then trivial if the problem is subject to Dirichlet-Dirichlet BCs (see
Table~\ref{tabR} and recall that $\cP_\ell^{(0,1)}(1)=1$ for each $\ell$).
On the other hand, if $\beta_b\neq 0$ then from (\ref{xinns}) and (\ref{tetanns}) with $(\alpha_a,\beta_a)=(1,0)$ we get
$$\xi_n = \nu_n(2\alpha_b + (n+2)^2\beta_b), \qquad \theta_n = -\nu_n (2\alpha_b + (n+1)^2\beta_b).$$
In addition, it is known that $\frac{d}{dx} \cP_n^{(0,1)} (1) = n(n+2)/2.$ With this information, one verifies 
that the polynomial at the right hand-side of (\ref{Rndm1}) satisfies the BC at $x=1$ and this completes the proof.\qed\\

\noindent The matrix $Q_N$ can be therefore written as
\begin{equation} \label{Qfact2}
Q_N = R_N^T \, \hat{F}_N R_N + \tilde{R}_N^T \, \tilde{G}_N \tilde{R}_N,
\end{equation}
see (\ref{R}) and (\ref{hFG}), where 
\begin{eqnarray} \label{hR}
\tilde{R}_N &=& \left(\begin{array}{cccccc}
       \xi_0\\
       \theta_0 & \xi_1 \\
       & \theta_1 & \ddots\\
       && \ddots & \xi_{N-1}\\
       &&& \theta_{N-1}\\
\end{array}\right) \in \Rr^{(N+1)\times N},\\
 \tilde{G}_N &=& \left(\tilde{g}_{mn}\right), \quad 
 \tilde{g}_{mn} = \int_{-1}^1 (1+x)^{2-\gamma} g(x) \cP_m^{(0,1)}(x) \cP_n^{(0,1)}(x) dx. 
 \label{tG}
 \end{eqnarray}
Now with an approach similar to the one considered in Proposition~\ref{propGFinf} and in the 
subsequent paragraph, 
which is essentially based on the recurrence relation for the  Jacobi polynomials $\cP_n^{(0,1)}$ \cite{Abram,EMOT}, 
we obtain that if we know the values of $\tilde{g}_{m0}$ for each $m=0,1,\ldots, 2N+1$ then
we can determine the remaining required values recursively. Moreover, 
if we let 
$$\tilde{\bfg}_0 = \left(\tilde{g}_{00},\tilde{g}_{10},\ldots\right)^T \in \ell_\infty$$ 
then we get, see \cite[16.4 formula (2)]{EMOT},
\begin{equation}\label{tg0}
\tilde{\bfg}_0 = g(\tilde{\mathcal{H}}) 
\left(\begin{array}{c}\tilde{g}_0\\ \tilde{g}_1\\\vdots \end{array}\right), \,\,\, 
\tilde{g}_m = \int_{-1}^{1} \frac{\cP^{(0,1)}_m(x)}{(1+x)^{\gamma-2}} dx =  \frac{(-1)^m 2^{3-\gamma} (\gamma-1)_m}{(3-\gamma)_{m+1}}, 
\end{equation}
where
\begin{eqnarray}\nonumber
\tilde{\mathcal{H}} &=& \left( \begin{array}{ccccc}
            \tilde{h}_{00} & \tilde{h}_{01}\\
            \tilde{h}_{10} & \tilde{h}_{11} & \tilde{h}_{12}\\
            &\tilde{h}_{21} & \tilde{h}_{22} & \tilde{h}_{23}\\
            && \ddots & \ddots & \ddots\\
            \end{array}\right),  \\ \nonumber
            \tilde{h}_{m,m-1} &=& \frac{m}{2m+1}, \quad \tilde{h}_{m,m} =\frac{1}{(2m+1)(2m+3)}, \quad \tilde{h}_{m,m+1} = \frac{m+2}{2m+3}.
\end{eqnarray} 

\begin{remark}\label{remg1}
If $\ga =1$ then $\tilde{g}_m=0$ for each $m\ge 1.$ Therefore, in this case,
$Q_N$ can be approximated up to machine precision with a banded matrix where the 
bandwidth depends on the number of numerically significant coefficients of the Legendre-Fourier series 
expansion of $g$ and $f.$ In particular, $Q_N$ is tridiagonal 
for the Boyd equation for which $g(x)\equiv g(-1)$ and $f\equiv 0.$ 
\end{remark}

\section{Error analysis and computation of corrected numerical eigenvalues.} \label{sec3}
In this section, we shall study the behavior of the error in the resulting numerical eigenvalues 
as $N$ increases and for a fixed index. \\
As usual, this is related to the regularity of the solution, namely, in this context,
to the regularity of the eigenfunctions. In particular, if $\gamma=0$ or 
$g(x) \equiv 0,$ then $q,$ and consequently $y,$  belongs to $C^\infty[-1,1].$
In this case, it is well-known that the errors in the approximations provided
by a spectral method decay exponentially. \\
Problems which require a deeper 
analysis are therefore those for which $q$ is unbounded at the left endpoint. 
We must observe that from (\ref{Qfact})--(\ref{hgmn}) and (\ref{Qfact2})--(\ref{tG})
one deduces that the spectral method we have derived is well defined for each $\ga \in (0,3),$ 
with $y(-1)=0$ if the left endpoint is singular. Nevertheless, we shall consider only   
the case $\ga \in (0,2]$  with $g(-1)\neq 0,$ namely only problems for which $x=-1$ is  a regular singular endpoint.
The generalization to essential singularities will be the topic of future research.\\
In this context, the results we are going to present are not only interesting from the theoretical point of view but 
they will also provide very simple, economical and effective 
techniques for the computation of corrected  numerical eigenvalues.\\

Let $\lambda^{(N)}$ be the approximation of the exact eigenvalue $\lambda$ as $N$ increases and let  
$y$ be the corresponding exact eigenfunction having the following expansion
\begin{eqnarray}
 y(x) = \sum_{n=0}^{+\infty} c_n \cR_n(x) &\equiv& \sum_{n=0}^{N-1} c_n \cR_n(x) +\sum_{n=N}^{+\infty} c_n \cR_n(x) \nonumber \\
 &\equiv& y_N(x) + \sum_{n=N}^{+\infty} c_n \cR_n(x). \label{espandiy}
\end{eqnarray}
By construction of $z_N$ and of $\lambda^{(N)},$ see (\ref{zN})-(\ref{wf}),  we can write
\begin{eqnarray*}
- \langle y_N,z_N''\rangle + \langle  y_N, q z_N\rangle = \lambda^{(N)} \langle y_N,z_N\rangle = 
  \lambda^{(N)} \left( \langle y,z_N\rangle + \langle y_N-y,z_N\rangle\right).
\end{eqnarray*}
On the other hand, it is evident that
\begin{eqnarray*}
- \langle z_N,y''\rangle + \langle  z_N, q y\rangle = \lambda \langle z_N,y\rangle.
\end{eqnarray*}
In addition
$$ \langle z_N,y''\rangle = \langle y,z_N''\rangle = \langle y_N,z_N''\rangle + \langle y-y_N,z_N''\rangle = \langle y_N,z_N''\rangle $$
since $z_N'' \in \Pi_{N-1}$ and $y-y_N \in \Pi_{N-1}^\bot.$
From these formulas we get 
\begin{eqnarray} \nonumber 
\lambda -\lambda^{(N)} &=& \frac{ \langle z_N, (q-\lambda^{(N)}) (y-y_N) \rangle}{\langle z_N,y\rangle} \\
\nonumber &=& \frac{ \langle z_N, q (y-y_N) \rangle}{\langle z_N,y\rangle} - \lambda^{(N)}\,\frac{\langle z_N, y-y_N\rangle}{\langle z_N,y\rangle}\\
&=& \frac{1}{\langle z_N,y\rangle} \left(\sum_{n=N}^{+\infty} c_n \langle \cR_n, q z_N\rangle\right) -\lambda^{(N)} \varepsilon_N,\label{error2}
\end{eqnarray}
being $\varepsilon_N = \langle z_N, y-y_N\rangle/\langle z_N,y\rangle.$
Therefore, an analysis of the behavior of the coefficients $c_n$  in (\ref{espandiy}) as $n$ increases is required and
the following result constitutes a first step.

\begin{proposition} If $n$ is sufficiently larger than the index of the eigenvalue, $\ga \in (0,2]$ and $g(-1)\neq 0$ then 
\begin{equation} \label{cnas}
c_n \approx -\frac{\langle \cR_n, (1+x)^{-\ga} gy\rangle}{a_{nn}} = - \frac{\int_{-1}^1 (1+x)^{-\ga} g(x) \cR_n(x) y(x) dx}{a_{nn}}. 
\end{equation}
\end{proposition}
\underline{Proof}: It is evident that 
$$a_{nn} c_n =  \langle \cR_n, (\lambda-q)  y\rangle.$$
We recall that we assumed $f$ in (\ref{vx}) to be analytical inside and over a   
Bernstein ellipse containing $[-1,1]$ and this implies that its Legendre coefficients 
decay exponentially. 
Thus, $\langle \cR_n,(\lambda -f)y\rangle$ becomes negligible with respect to $\langle \cR_n, (1+x)^{-\ga} g y\rangle$
as $n$ increases and this complete the proof. \qed\\

From (\ref{error2}), by using similar arguments,  one deduces that the main contribution to the error in the eigenvalue 
approximation is given by
\begin{eqnarray}\label{error3_0}
\lambda-\lambda_N &\approx& \frac{1}{\langle z_N,y\rangle} \sum_{n=N}^{+\infty} c_n \langle \cR_n, (1+x)^{-\ga} g z_N\rangle\\
 \label{error3} &\approx& -\frac{1}{\langle z_N,y\rangle} \sum_{n=N}^{+\infty} \frac{\left\langle\cR_n, (1+x)^{-\ga} g z_N\right\rangle\, \langle \cR_n, (1+x)^{-\ga} g y\rangle}{a_{nn}}.
\end{eqnarray}

We recall the following asymptotic estimate \cite{Sid09}.
\begin{proposition}\label{propsi} Let $\psi \in C^{\infty}(-1,1)\bigcap C[-1,1]$ have the expansion
$$ \psi(x) = \psi(-1) \sum_{j=0}^L \psi_j (1+x)^{\sigma_j} + O((1+x)), \qquad \mbox{as \,\,} x \rightarrow -1^+,$$
with $\psi(-1)\neq 0, \psi_0=1$ and $\sigma_0=0< \sigma_1 < \ldots< \sigma_L < 1.$ If $s\in (-1,+\infty)\setminus \Nn_0$ 
and if $n$ is sufficiently large then
\begin{eqnarray} \label{asint}
\lefteqn{\int_{-1}^{1} (1+x)^s \psi(x) \cP_n(x) dx =}\\
&&\psi(-1) \left(\sum_{j=0}^L \psi_j \int_{-1}^{1} (1+x)^{s+\sigma_j} \cP_n(x) dx\right) + O\left(n^{-2s-4}\right). \qquad \mbox{\qed} \label{asint1}
\end{eqnarray}
\end{proposition}

We can now prove the following result which concerns (weakly) regular problems not subject to the Dirichlet boundary condition at the left endpoint. 
\begin{theorem}\label{tecnod}
If $y(-1)g(-1)\neq 0,$ $\ga \in (0,1)$ and if $N$ is sufficiently larger than the 
index of the eigenvalue then 
\begin{equation}\label{errnodiri}
\lambda-\lambda^{(N)} \approx - \,\frac{\omega ^2 g^2(-1) z_N(-1) y(-1) }{\langle z_N,y \rangle \,p\, (N+1)^p} 
\end{equation}
where 
\begin{eqnarray}\label{ordnodiri}
\omega &=& \frac{2^{2-\gamma} \Gamma(3-\gamma)}{(1-\gamma)\Gamma(\gamma)}, \qquad p = 6-4\ga. 
\end{eqnarray}
\end{theorem}
\underline{Proof}:
From (\ref{Rnp2}), (\ref{error3}) and the previous Proposition with $s=-\gamma,$ $\psi = g\,z_N$ or  $\psi=g y,$ and $L=0$ we obtain
$$ \lambda-\lambda^{(N)} \approx -\frac{g^2(-1) z_N(-1) y(-1) }{\langle z_N,y \rangle} \left(\sum_{n=N}^{+\infty} \frac{\langle \cR_n, (1+x)^{-\gamma}\rangle^2}{a_{nn}}\right).$$
Now, we can rewrite the last equation in (\ref{hg0}) as follows 
$$ \langle \cP_n, (1+x)^{-\gamma}\rangle = \frac{(-1)^n 2^{1-\ga} (\ga)_n}{(1-\ga)_{n+1}}  = 
   \frac{(-1)^n\,2^{1-\ga}\, \Gamma(1-\ga)\,\Gamma(\ga + n)}{\Gamma(\ga)\,\Gamma(2-\ga+n)}.$$
Therefore, by using the following expansion of the ratio of two gamma functions 
\begin{equation}\label{ratgamma}
\frac{\GG(z+a)}{\GG(z+b)} = z^{a-b} \left(1 + \frac{(a-b)(a+b-1)}{2z} + O (z^{-2})\right),\quad
z\gg 0,
\end{equation}
with $z=n+1/2,$ we get
$$ \langle \cP_n, (1+x)^{-\gamma}\rangle = 
\frac{(-1)^n 2^{1-\ga} \Gamma(1-\ga)}{\Gamma(\ga)} \left(n+\frac{1}{2}\right)^{2\ga-2} \left(1 + O\left(\frac{1}{n^2}\right)\right)$$
and, consequently,  
\begin{eqnarray*}
\lefteqn{\frac{(-1)^n\Gamma(\ga)}{2^{1-\ga} \Gamma(1-\ga)} \,\langle \cR_n, (1+x)^{-\ga}\rangle }\\
&\approx&\left(\xi_n \left(n+\frac{1}{2}\right)^{2\ga-2} - \eta_n \left(n+\frac{3}{2}\right)^{2\ga-2} + \theta_n \left(n+\frac{5}{2}\right)^{2\ga-2}\right) \\ 
&\approx&\left(n+\frac{3}{2}\right)^{2\ga-2} \left(\xi_n-\eta_n+\theta_n  - (2\ga-2)(\xi_n-\theta_n)\left(n+\frac{3}{2}\right)^{-1}\right).
\end{eqnarray*}
We recall that if $n$ is sufficiently large then $\xi_n = 1,$ see (\ref{xias}). In addition,
by using (\ref{xins})--(\ref{etas}) it is possible to verify with some computations that if $\beta_a \neq 0$ then
$$ \xi_n-\eta_n+\theta_n = \frac{8}{2n+3} \left(1 + O\left(\frac{1}{n}\right)\right), \qquad \xi_n -\theta_n = 2 \left(1 + O\left(\frac{1}{n}\right)\right).$$
Hence, see (\ref{annas}) and (\ref{ordnodiri}),
\begin{eqnarray}
\nonumber \langle \cR_n, (1+x)^{-\ga}\rangle  &\approx& \frac{(-1)^n 2^{1-\ga} \Gamma(1-\ga) (8 -4\ga)}{\Gamma(\ga)} \left(n+\frac{3}{2}\right)^{2\ga-3}\\
                                    &=& (-1)^n 2\omega \left(n+\frac{3}{2}\right)^{-p/2} \nonumber,\\
\label{rnmg} a_{nn}^{-1} \langle \cR_n, (1+x)^{-\ga}\rangle &\approx& (-1)^n\,\frac{\omega}{2} \left(n+\frac{3}{2}\right)^{-p/2-1}. 
\end{eqnarray}
Therefore
\begin{eqnarray*}
\lambda-\lambda^{(N)} &\approx& -\frac{g^2(-1) z_N(-1) y(-1)}{\langle z_N,y \rangle} \left(\sum_{n=N}^{+\infty} \frac{\langle \cR_n, (1+x)^{-\gamma}\rangle^2}{a_{nn}}\right) \\
&\approx&   -\,\frac{\omega^2\,g^2(-1) z_N(-1) y(-1)}{\langle z_N,y \rangle} \sum_{n=N}^{+\infty} \left(n+3/2\right)^{-p-1} \\
&\approx& -\,\frac{\omega^2\,g^2(-1) z_N(-1) y(-1)}{\langle z_N,y \rangle} \int_{N}^{+\infty} \left(n+1\right)^{-p-1} dn\\
&=& - \,\frac{\omega^2\,g^2(-1) z_N(-1) y(-1)}{\langle z_N,y \rangle \,p\,(N+1)^{p}}
\end{eqnarray*}
which is the statement of the theorem. \qed\\

This result immediately suggests a very simple formula for the correction of the numerical eigenvalue. First of 
all, we assume the numerical and the exact eigenfunctions have been normalized so that  
\begin{eqnarray}\label{normal}
\langle z_N,z_N\rangle = \bfzeta_N^T B_N \bfzeta_N &=&1, \qquad \langle y,y \rangle = 1,\\
z_N(-1)&>&0, \qquad \,\,y(-1)>0. \nonumber
\end{eqnarray}
By using the orthogonality of the Legendre polynomials, this permits the estimates
$y(-1) \approx z_N(-1),$ $\langle z_N,y \rangle \approx 1,$
and consequently, see (\ref{errnodiri}), 
$$ \lambda-\lambda^{(N)} \approx - \frac{\left(\omega \,g(-1)\, z_N(-1) \right)^2 }{p\, (N+1)^p}.$$ 
In addition, we observe that
the term $\lambda^{(N)} \varepsilon_N$ in (\ref{error2}) can be of some relevance if $N$ is not so much larger than the 
index of the eigenvalue. By virtue of (\ref{cnas}), (\ref{asint})-(\ref{asint1}) and of (\ref{rnmg}), we therefore decided to consider 
the following approximation of $\varepsilon_N$ which can be computed with a very low cost
\begin{eqnarray}
\nonumber \varepsilon_N &\approx & \langle z_N,y-y_N\rangle 
= c_N \langle z_N, \cR_N \rangle + c_{N+1} \langle z_N,\cR_{N+1}\rangle \\
&\approx&  \bar{c}_N \langle z_N, \cR_N \rangle + \bar{c}_{N+1} \langle z_N,\cR_{N+1}\rangle
       \equiv \bar{\varepsilon}_N, \label{vareps}
\end{eqnarray}
where 
\begin{equation}\label{cnbarnd}
\bar{c}_n \equiv - \frac{(-1)^n \omega\,g(-1)\,z_N(-1)}{2}  \left(n +\frac{3}2\right)^{-p/2-1}.
\end{equation}
It is worth recalling that 
$\langle z_N,\cR_N \rangle = b_{N,N-2}\zeta_{N-2,N} + b_{N,N-1}\zeta_{N-1,N},$ and
$\langle z_N,\cR_{N+1} \rangle = b_{N+1,N-1}\zeta_{N-1,N}$ (see (\ref{geneig}),(\ref{abmn}) and (\ref{Bfact})).\\

All these arguments lead to the following formula for the correction of the numerical eigenvalue to be used
when the BC at the left endpoint is not of Dirichlet type
\begin{equation}\label{lamcnod}
\mu^{(N)} \equiv \lambda^{(N)} (1-\bar{\varepsilon}_N)- \frac{(\omega\, g(-1)\,z_N(-1))^2}{p \,(N+1)^p}   \approx \lambda.
\end{equation}

The main steps of the procedure for the computation of such $\lambda^{(N)}$ and $\mu^{(N)}$  are summarized in Algorithm~\ref{algnod}. 
With respect to the notation we have used so far, we add the further index $k$ which represents the index of the eigenvalue. 
Its value belongs to $\left\{1,\ldots,M\right\}$ being $M$ the number of smallest eigenvalues requested.\\

\begin{algorithm}
\caption{Solution of a (weakly) regular problem with $y(-1)\neq 0.$} \label{algnod}
\textbf{Input:} $f,g,\ga,(\alpha_a,\beta_a),(\alpha_b,\beta_b),M,N$\smallskip\\ 
\textbf{Require:} $\ga \in (0,1),\,\beta_a \neq 0$ and $M\leq N$ \smallskip\\ 
\textbf{Output:} $\lambda_k^{(N)}$ and $\mu_k^{(N)}$ \,\,\mbox{for } $k=1,\ldots,M$ \smallskip 
\begin{algorithmic}[1]
\State Compute $A_N,$ $B_N$ and $Q_N$ by using (\ref{amm1}),(\ref{Bfact})--(\ref{hFG})
\State Solve the generalized eigenvalue problem 
$$
(A_N+Q_N) \bfzeta_N^{(k)} = \lambda_{k}^{(N)} B_N \bfzeta_N^{(k)},  \qquad k=1,\ldots,M,
$$
with $\left(\bfzeta_N^{(k)}\right)^T B_N \bfzeta_N^{(k)} = 1;$
\For{$k \gets 1,M$}
 \State  Compute $z_{k,N}(-1) = \sum_{n=0}^{N-1} \zeta_{n,N}^{(k)} \cR_n(-1)$
 \State Use (\ref{ordnodiri}),(\ref{vareps})--(\ref{lamcnod}) to determine $\mu_k^{(N)}.$
\EndFor
\end{algorithmic}

\end{algorithm}

Let us now consider problems subject to the Dirichlet boundary condition at the left endpoint with $\ga \in (0,2).$ 
In this case, Proposition~\ref{propsi} with $s=-\ga$
is not directly applicable to the inner products in (\ref{error3}). In fact, we need $s> -1$  and $\psi(-1)\neq 0.$ 
Nevertheless, it is evident that since $z_N(-1)=y(-1)=0$  we can write
$$z_N(x) = (1+x) \hat{z}_N (x), \qquad y(x) = (1+x) \hat{y}(x),$$
with $\hat{z}_N(-1) = z_N'(-1)\neq 0$ and $\hat{y}(-1) = y'(-1)\neq 0.$
Consequently
\begin{eqnarray*}
\langle \cR_n, (1+x)^{-\ga} g z_N \rangle &=& \langle \cR_n, (1+x)^{1-\ga} g \hat{z}_N \rangle, \\
\langle \cR_n, (1+x)^{-\ga} g y \rangle &=& \langle \cR_n, (1+x)^{1-\ga} g \hat{y} \rangle.
\end{eqnarray*}
Now, an analysis of the behavior of $\hat{y}$ in proximity of $x=-1$ is required 
for the application of (\ref{asint})-(\ref{asint1}) and 
a Frobenius-type method provides the following result. 
\begin{lemma} If $y(-1) = 0,$ $g(-1)\neq 0$ and $\ga \in (0,2)$ then an eigenfunction $y$ admits 
the following expansion as $x\rightarrow -1^+:$
\begin{equation} \label{frobe}
y(x) = (1+x) y'(-1) \left[\sum_{j=0}^L \chi_j (1+x)^{j(2-\ga)}  + O((1+x)^{s_\ga})\right]. 
\end{equation}
Here
\begin{eqnarray}
 L &=& \max\left(0,\lceil (\ga-1)/(2-\ga)\rceil\right), \qquad s_\ga \ge 1, \label{L}\\
\chi_0 &=& 1, \qquad \chi_{j+1} = \frac{g(-1)}{(j+1)(2-\ga)(1+(j+1)(2-\ga))} \chi_j, \label{chi} 
\end{eqnarray}
i.e.
$$ \chi_j = \frac{1}{\left(\frac{3-\ga}{2-\ga}\right)_j \, j!} \left(\frac{g(-1)}{(2-\ga)^2}\right)^j.$$
\end{lemma}

\begin{remark} 
It must be observed that the term with the summation in (\ref{frobe}) represents the truncation of a fractional power series expansion
at $x=-1$ of the solution of (\ref{slp}) with 
\begin{equation}\label{pref}
\lambda=0, \quad q(x) = g(-1)/(1+x)^\ga, \quad y(-1) = 0, \quad y'(-1) \quad \mbox{assigned}.
\end{equation}
In fact, if one sets $t=(1+x)^{1-\ga/2}$ and $y(x)= u(t)$ then one gets
$$ u'' - \frac{(2\nu-1)}t u' - 4 \nu^2 g(-1) u = 0, \qquad \nu = 1/(2-\ga).$$
A solution of this equation subject to $u(0)=0$ is proportional to \cite{Bow58}
\begin{eqnarray*}
t^{\nu} J_\nu\left(2 {\rm i}\nu \sqrt{g(-1)} \right) &=& 
\frac{t^{2\nu}\left({\rm i} \nu \sqrt{g(-1)}\right)^\nu}{\Gamma(\nu+1)} \,_0F_1\left(\nu+1; \nu^2 g(-1) t^2\right)\\
&\propto& t^{2/(2-\ga)} \,_0F_1\left(\frac{3-\ga}{2-\ga}; \frac{g(-1) t^2}{(2-\ga)^2}\right).
\end{eqnarray*}
Here $J_\nu$ is the Bessel function of the first kind and $\,_0F_1$ a confluent hypergeometric limit function. 
Therefore, the solution of the initial value problem (\ref{slp})-(\ref{pref}) is
\begin{eqnarray*}
y(x) &=& (1+x) y'(-1) \,_0F_1\left(\frac{3-\ga}{2-\ga}; \frac{g(-1) (1+x)^{2-\ga}}{(2-\ga)^2}\right)\\
     &=& (1+x) y'(-1) \sum_{j=0}^{+\infty} \chi_j (1+x)^{j(2-\ga)}.
\end{eqnarray*}
\end{remark}

\begin{remark}
For the special value $\ga=1,$ like the Boyd equation, it is possible to verify that a solution of (\ref{slp}) subject to $y(-1)=0$
admits a (classical) power series expansion at $x=-1$. Moreover, the coefficients
of the expansion of $y$ in (\ref{espandiy}) decays exponentially and, see Remark~\ref{remg1}, 
the matrix $Q_N$ can be approximated up to machine 
precision with a banded one with bandwidth independent of $N.$  From these observations, 
we deduce that if $\ga =1$ then the errors in the approximation of the eigenvalues decay exponentially with respect to $N.$
\end{remark}
Problems to be studied are therefore those with $\ga \in (0,2)\setminus\left\{1\right\}$ and this is done in the next theorem.

\begin{theorem}
If $y(-1)=0,$ $g(-1) \neq 0$ and $\ga \in (0,2)\setminus\left\{1\right\}$ then 
\begin{eqnarray}\nonumber
\lambda-\lambda^{(N)} 
   &\approx& - \frac{g^2(-1) z_N'(-1) y'(-1) }{\langle y,z_N \rangle(N+1)^{p}} \, \sum_{j=0}^L \frac{\omega_j}{(p+2j(2-\ga)) (N+1)^{2j(2-\ga)}} \label{errdiri}
\end{eqnarray}
where $L$ is defined in (\ref{L}) and, see (\ref{chi}), 
\begin{eqnarray}\label{orddiri}
p &=& 10-4\ga, \\
\label{omega}
\omega_j &=& \frac{2^{4-\ga} \, \Gamma (3-\ga)}{\Gamma(\ga-1)} \,\hat{\omega}_j, \,\,  
\hat{\omega}_j = \frac{2^{(2-\ga)(j+1)}\,\Gamma(1 + (2-\ga)(j+1))} {\Gamma(1-(2-\ga)(j+1))} \,\chi_j.
\end{eqnarray}
\end{theorem}
\underline{Proof}: 
From the premises of this theorem, Propositions~\ref{propU} and \ref{propsi}, (\ref{Rndm1}), (\ref{tg0}), (\ref{ratgamma})
 with $z=n+1$, (\ref{orddiri}),
and (\ref{xias})-(\ref{etas})
we obtain that if $n\geq N$ 
with $N$ sufficiently large then
\begin{eqnarray*}
\lefteqn{\llangle \cR_n, (1+x)^{-\ga} g z_N \rrangle = \llangle \cR_n, (1+x)^{1-\ga} g \hat{z}_N \rrangle } \\
&\approx& g(-1) \hat{z}_N(-1) \llangle \cR_n, (1+x)^{1-\ga} \rrangle = g(-1) z_N'(-1) \llangle \cU_n, (1+x)^{2-\ga} \rrangle\\
&\approx& \frac{(-1)^n  g(-1) z_N'(-1) 2^{3-\ga}\Gamma(3-\ga)}{\Gamma(\ga-1)} \left(\frac{\xi_n}{(n+1)^{p/2}} - \frac{\theta_n}{(n+2)^{p/2}}\right)\\
&\approx& \frac{(-1)^n  g(-1) z_N'(-1) 2^{4-\ga}\Gamma(3-\ga)}{\Gamma(\ga-1)} \left(n + \frac{3}2\right)^{-p/2}.
\end{eqnarray*}
Similarly, by using also (\ref{omega}) and the previous lemma, we get
$$\llangle \cR_n, (1+x)^{1-\ga} g \hat{y} \rrangle
\approx g(-1) y'(-1) \sum_{j=0}^L \chi_j \,\llangle \cU_n, (1+x)^{(2-\ga)(j+1)} \rrangle
$$
with 
\begin{eqnarray*}
\chi_j \left\langle \cU_n, (1+x)^{(2-\ga)(j+1)} \right\rangle 
&\approx& (-1)^n 4 \hat{\omega}_j \left(n + \frac{3}2\right)^{-p/2-2j(2-\ga)}.
\end{eqnarray*}
Hence, by considering (\ref{annas}) and (\ref{cnas}) we obtain
\begin{equation}\label{cnasd}
c_n \approx - (-1)^ng(-1)y'(-1) \sum_{j=0}^L \hat{\omega}_j \left(n + \frac{3}2\right)^{-p/2-1-2j(2-\ga)} \equiv \bar{c}_n
\end{equation}
so that, see (\ref{error2}),
$$\lambda -\lambda_N \approx - \frac{g^2(-1) z_N'(-1) y'(-1)}{\langle z_N,y\rangle} \sum_{j=0}^L\sum_{n=N}^{+\infty} \omega_j \left(n+\frac{3}2\right)^{-p-1-2j(2-\ga)}.$$
The statement follows by using an integral estimate. \qed\\

Let us now discuss how one can use this result for the correction of the numerical eigenvalues. 
Following the idea used for problems with $y(-1)\neq0,$ we consider the normalization specified in (\ref{normal}),
with $z'_N(-1)>0 $ and $y'(-1)>0,$
by which we get $\langle z_N,y \rangle \approx 1.$ On the other hand, the estimate $y'(-1)\approx z_N'(-1)$
turns out to be rather poor. In fact, we have just established that if $n$ is sufficiently large 
then $c_n = O\left((n+3/2)^{-p/2-1}\right)$ and it is possible to verify that 
$$\cR_n'(-1) = \cU_n(-1)= (-1)^n (2n+3)+ O(1).$$ 
Therefore, $ c_n \cR_n'(-1) = O\left((n+3/2)^{-p/2}\right)$ approaches zero rather slowly when $p \rightarrow 2^+,$ i.e. $\ga \rightarrow 2^-$. 
For example, if $\ga = 1.9$ then $p/2 = 6/5.$ By considering (\ref{zN}) and (\ref{espandiy}), the following approximation 
$$ y'(-1) = y_N'(-1) + \sum_{n=N}^\infty c_n \cR_n'(-1) \approx z_N'(-1) + \sum_{n=N}^\infty c_n \cR_n'(-1)$$
turns out to be more appropriate. Now, from (\ref{cnasd}) we obtain
\begin{eqnarray}
\nonumber \sum_{n=N}^{+\infty} c_n \cR_n'(-1) &\approx& - g(-1) y'(-1) \sum_{n=N}^{+\infty} \sum_{j=0}^L 2\hat{\omega}_j \left(n+\frac{3}2 \right)^{-1-2(2-\ga)(j+1)}\\
\label{dN0} &\approx & -\left(\frac{g(-1)}{2-\ga} \sum_{j=0}^L \frac{\hat{\omega}_j}{(j+1)(N+1)^{2(2-\ga)(j+1)}}\right) y'(-1)\\
\label{dN} &\equiv& - d_N\,y'(-1),
\end{eqnarray}
and, consequently, $y'(-1) \approx z_N'(-1) - d_N y'(-1),$ i.e. $y'(-1) \approx z_N'(-1)/(1+d_N)$.\\

Summarizing, if the problem is subject to the Dirichlet BC at the singular endpoint and if $\ga \in(0,2)\setminus\left\{1\right\}$ then 
we correct the numerical eigenvalues as follows
\begin{equation}\label{lamcdiri}
\mu^{(N)} \equiv \lambda^{(N)} (1-\bar{\varepsilon}_N)- 
\frac{\left(g(-1) z_N'(-1)\right)^2}{(1+d_N) \,(N+1)^{p}}  \sum_{j=0}^L \frac{\omega_j}{(p+2j(2-\ga)) (N+1)^{2j(2-\ga)}}
\end{equation}
where  $\bar{\varepsilon}_N$ is defined in (\ref{vareps}) with $\bar{c}_N$ and $\bar{c}_{N+1}$ 
given by (\ref{cnasd}) with $n=N,N+1,$ respectively. The main steps of the procedure for their computation are 
listed in Algorithm~\ref{algdiri}.\\

\begin{algorithm}
\caption{Solution of a problem with $y(-1)=0$ and $\ga \in(0,2)\setminus\left\{1\right\}.$} \label{algdiri}
\textbf{Input:} $f,g,\ga,(\alpha_b,\beta_b),M,N$ \smallskip\\
\textbf{Require:} $\ga \in(0,2)\setminus\left\{1\right\}$  and $M\leq N$ \smallskip \\
\textbf{Output:} $\lambda_k^{(N)}$ and $\mu_k^{(N)}$ \,\,\mbox{for } $k=1,\ldots,M$ 

\begin{algorithmic}[1]
%
%
\State Set $(\alpha_a,\beta_a)=(1,0)$
\State Compute $A_N,$ $B_N$ and $Q_N$ by using (\ref{amm1}),(\ref{Bfact}) and (\ref{Qfact2})
\State Solve the generalized eigenvalue problem $$(A_N+Q_N) \bfzeta_N^{(k)} = \lambda_{k}^{(N)} B_N \bfzeta_N^{(k)}, \qquad k=1,\ldots,M,$$
with $\left(\bfzeta_N^{(k)}\right)^T B_N \bfzeta_N^{(k)} = 1$
\State Set $L = \mbox{max}(0,\lceil(\ga-1)/(2-\ga)\rceil)$
\State Use (\ref{chi})-(\ref{omega}) to compute $\hat{\omega}_j$ and $\omega_j$  for $j=0,\ldots,L$
\State Determine $d_N$ defined in (\ref{dN0})-(\ref{dN})
\For{$k \gets 1,M$}
 \State Compute $z_{k,N}'(-1) = \sum_{n=0}^{N-1} \zeta_{n,N}^{(k)} \cR'_n(-1)$
 \State Use (\ref{orddiri})-(\ref{omega}) and (\ref{lamcdiri}) to determine $\mu_k^{(N)}$
\EndFor
\end{algorithmic}
\end{algorithm}

The final error analysis we are going to present concerns problems with $\ga=2$ and $g(-1)>0.$  In this case, the application of the 
Frobenius method allows to state that the exact eigenfunction satisfies 
\begin{equation}\label{frobeg2}
y(x) \equiv (1+x)\hat{y}(x) = \chi\, (1+x)^\varrho (1 + O(1+x)), \qquad \mbox{as } x \rightarrow -1^+, 
\end{equation}
where $\chi$ is a free parameter while $\varrho$ is the positive root of the indicial equation $\varrho^2 -\varrho -g(-1),$ i.e.
\begin{equation}\label{ro}
 \varrho = \frac{1 + \sqrt{1+4g(-1)}}{2} > 1.
\end{equation}
For instance, if $g(x)\equiv g(-1)$ and $f(x)\equiv 0$ then a solution of (\ref{slp}) subject to $y(-1)=0$ is 
proportional to \cite{Bow58}
$$ \sqrt{1+x}\, J_{\varrho-0.5} \left(\sqrt{\lambda}\,(1+x)\right) \propto (1+x)^\varrho \,_0F_1\left(\varrho +\frac{1}{2}; -\frac{\lambda}{4}\,(1+x)^2\right).$$
We shall proceed by assuming the exact eigenfunction has been normalized so that 
$$ 
 \langle y,y\rangle = 1 \qquad \mbox{and} \qquad \chi > 0.
$$
Concerning the corresponding numerical eigenfunction, we assume 
\begin{equation}\label{normalzN}
 \langle z_N,z_N\rangle = 1 \qquad \mbox{and}  \qquad \lim_{N\rightarrow +\infty} \frac{\hat{z}_N(-1)}{\hat{y}_N(-1)} = \kappa > 0
\end{equation}
being 
\begin{equation}\label{hzhyn}
\hat{z}_N(x) = \sum_{n=0}^{N-1}\zeta_{n,N}\,\cU_n(x), \qquad \hat{y}_N(x) = \sum_{n=0}^{N-1}c_n\,\cU_n(x).
\end{equation}
In particular, the second formula in (\ref{normalzN}) state that $\hat{z}_N(-1)$ and $\hat{y}_N(-1)$ 
are infinitesimal of the same order 
as $N$ increases.\\ 
With these notations, we can prove the following result.

\begin{theorem}\label{teoconvg2}
If $\ga = 2,$ $y(-1)=0,$ $\varrho \notin \Nn,$ see (\ref{ro}), and if $N$ is sufficiently larger then the index of the 
eigenvalue then
\begin{eqnarray} \label{errg2}
 \lambda-\lambda^{(N)} &\approx& \frac{2^{\varrho} \Gamma(\varrho+1) (\varrho -1) \chi\, \hat{z}_N(-1)}{\Gamma(1-\varrho)} (N+1)^{-2\varrho}\\
                       &\approx& -\frac{2}{\kappa} \left(\frac{(\varrho-1)\,\hat{z}_N(-1)}{N+1}\right)^2  = O((N+1)^{-p}) \label{errg2_1}
\end{eqnarray}
where $\kappa$ is the limit value in (\ref{normalzN}) and 
\begin{equation}\label{ordg2}
 p = 4\varrho-2 =2\sqrt{1 + 4g(-1)}.
\end{equation}
\end{theorem}
\underline{Proof}: We begin by studying the asymptotic behavior of the coefficients $c_n$ with 
an approach similar to the one used in the proof of the previous theorems.
From (\ref{cnas}), (\ref{asint})-(\ref{asint1}) and (\ref{frobeg2}) we get 
\begin{eqnarray*}
c_n &\approx& - \frac{\chi\,g(-1) \langle \cR_n, (1+x)^{\varrho-2}\rangle}{a_{nn}} = - \frac{\chi\,g(-1) \langle \cU_n, (1+x)^{\varrho-1}\rangle}{a_{nn}}\\
    &=& - \frac{(-1)^n \chi\,g(-1) 2^{\varrho+1} \Gamma(\varrho)}{a_{nn} \Gamma(2-\varrho)}\left(n+\frac{3}{2}\right)^{1-2\varrho} \left(1+O(n^{-2})\right)\\
    &\approx&  \frac{(-1)^n \chi 2^{\varrho-1} \Gamma(\varrho+1)}{\Gamma(1-\varrho)}\left(n+\frac{3}{2}\right)^{-2\varrho}.
\end{eqnarray*}
In particular, for the last estimate we used the fact that $g(-1) = \varrho(\varrho-1)$ and (\ref{annas}). \\

Let us now consider $\langle \cR_n, (1+x)^{-2} g z_N\rangle=\langle \cU_n, g  \hat{z}_N\rangle.$ 
If we let $g(x) = g(-1) + (1+x) \tilde{g}(x),$ with $\tilde{g} \in C^{\infty}([-1,1]),$ 
then we obtain
\begin{eqnarray*}
\langle \cU_n, g  \hat{z}_N\rangle
&=& g(-1) \langle \cU_n, \hat{z}_N\rangle + \langle \cU_n, (1+x) \tilde{g} \hat{z}_N\rangle\\
&\approx& \varrho(\varrho-1) \langle \cU_n, \hat{z}_N\rangle
\end{eqnarray*}
with a remainder that decreases exponentially. In addition, if one uses  \cite[16.4, formula(20)]{EMOT} 
then one gets
$ \langle \cP_{m}^{(0,1)},\cP_{n}^{(0,1)} \rangle = P_m^{(0,1)} (-1) \int_{-1}^1 \cP_n^{(0,1)} (x) dx= (-1)^n\,2/(n+1)$ for each  $m \leq n.$
Clearly, this implies $ \langle \cU_{m},\cU_{n} \rangle = \cU_m (-1) \int_{-1}^1 \cU_n (x) dx,$ \mbox{for each } $m < n$
and consequently 
\begin{eqnarray*}
\langle \cR_n, (1+x)^{-2} g z_N\rangle &\approx& \varrho(\varrho-1) \hat{z}_N(-1) \int_{-1}^1 \cU_n(x) dx\\
&\approx& \frac{ (-1)^n 8 \varrho(\varrho-1) \hat{z}_N(-1) }{2n+3}, \qquad \forall n\ge N.
\end{eqnarray*}
Therefore, the first estimate (\ref{errg2}) follows from (\ref{error3_0}),  with $\langle z_N,y\rangle \approx 1,$ 
and from the application of an integral estimate. It must be said that this result is only a starting point 
and it is not particularly useful both from the theoretical point of view and for the derivation of a correction formula. 
This is because $\hat{z}_N(-1)$ approaches zero as $N$ increases
and if we don't know its infinitesimal order then we don't know the order of convergence of $\lambda-\lambda^{(N)}.$ Concerning the 
computation of a corrected numerical eigenvalue we need an estimate of $\chi$. Let us therefore discuss the approximations in 
(\ref{errg2_1}). By using the formula obtained for $c_n$ and by observing that $\cU_n(-1) = (-1)^n (2n+3) + O(1),$ one gets
\begin{eqnarray*}
0= y'(-1) &=& \hat{y}_N(-1) + \sum_{n=N}^{+\infty} c_n \cU_n(-1)\\
&\approx& \hat{y}_N(-1) + \frac{\chi\, 2^{\varrho}\,\Gamma(\varrho+1)}{\Gamma(1-\varrho)} \sum_{n=N}^{+\infty} \left(n+\frac{3}{2}\right)^{1-2\varrho}\\
&\approx& \frac{\hat{z}_N(-1)}{\kappa} + \frac{\chi\, 2^{\varrho-1}\,\Gamma(\varrho+1)}{\Gamma(1-\varrho) (\varrho-1)} (N+1)^{2-2\varrho}
\end{eqnarray*}
where $\kappa$ is defined in (\ref{normalzN}). This implies $\hat{z}_N(-1) = O((N+1)^{2-2\varrho})$ and consequently 
$$\lambda-\lambda^{(N)} = O((N+1)^{-p}), \qquad p = 4\varrho-2.$$
Moreover
$$\chi \approx -\frac{\Gamma(1-\varrho) (\varrho-1) \hat{z}_N(-1)}{\kappa\,\Gamma(\varrho +1) 2^{\varrho-1}} (N+1)^{2\varrho-2}$$
which, with a simple substitution, completes the proof of (\ref{errg2_1}). \qed\\

We must now spend few words concerning the limit value $\kappa.$ It is evident that intuitively one would expect 
$\kappa=1,$ namely $\hat{z}_N(-1) \approx \hat{y}_N(-1).$ Nevertheless, the results of several numerical experiments
we have conducted by considering different $g,$ $f$ and different boundary conditions at  $x=1$
indicate that this is not the case. In particular, such tests lead us to the following assumption
\begin{equation}\label{kapcong}
\kappa = \lim_{N\rightarrow +\infty}  \frac{\hat{z}_N(-1)}{\hat{y}_N(-1)} = (2\varrho-1)/\varrho^2.
\end{equation}
Currently, this is an experimental result but 
we didn't a find a problem for which it doesn't work. 
By using it, we get that we can correct the numerical eigenvalues with the following very simple formula
\begin{equation}\label{mug2}
\mu^{(N)} = (1-\bar{\varepsilon}_N)\lambda^{(N)} - \frac{2}{2\varrho-1} \left( \frac{\varrho(\varrho-1) \hat{z}_N(-1)}{N+1}\right)^2,
\end{equation}
where, once again, $\bar{\varepsilon}_N$ is defined in (\ref{vareps}) with 
\begin{equation}\label{CNg2}
\bar{c}_M     = (-1)^{M+1} \frac{(\varrho-1)\varrho^2 \hat{z}_N(-1)}{(2\varrho-1) (N+1)^2} \left(\frac{2N+2}{2M+3}\right)^{2\varrho}, \qquad M=N,N+1.
\end{equation}
The procedure for their computation is sketched in Algorithm~\ref{algg2}.

\begin{algorithm}
\caption{Solution of a problem with $y(-1)=0$ and $\ga =2.$} \label{algg2}
\textbf{Input}: $f,g,(\alpha_b,\beta_b),M,N$\smallskip\\
\textbf{Require}: $g(-1) >0,$ and $M\leq N$ \smallskip\\
\textbf{Output}: $\lambda_k^{(N)}$ and $\mu_k^{(N)}$ \,\,\mbox{for } $k=1,\ldots,M$
\begin{algorithmic}[1]
%
%
\State Set $\ga =2$ and $(\alpha_a,\beta_a)=(1,0)$
\State Compute $A_N,$ $B_N$ and $Q_N$ by using (\ref{amm1}),(\ref{Bfact}) and (\ref{Qfact2})
\State Solve the generalized eigenvalue problem $$(A_N+Q_N) \bfzeta_N^{(k)} = \lambda_{k}^{(N)} B_N \bfzeta_N^{(k)}, \qquad k=1,\ldots,M,$$
with $\left(\bfzeta_N^{(k)}\right)^T B_N \bfzeta_N^{(k)} = 1$
\State Set $\varrho = \left(1+\sqrt{1+4g(-1)}\right)/2$
\For{$k \gets 1,M$}
 \State Compute $\hat{z}_{k,N}(-1)=\sum_{n=0}^{N-1} \zeta_{n,N}^{(k)} \cU_n(-1)$
 \State Use (\ref{mug2})-(\ref{CNg2}) and (\ref{vareps}) to determine $\mu_k^{(N)}$
\EndFor
\end{algorithmic}
\end{algorithm}

\section{Numerical tests} \label{sec4}
The method described and the algorithms for the a posteriori correction 
were implemented in Matlab ({\tt ver.R2017a}). In particular, we solved the arising 
generalized eigenvalue problem (\ref{geneig}) 
by using the {\tt eigs} function with option ``{\tt SM}'' for getting the ones of smallest magnitude.
In addition, as we already said in Section~\ref{secQ}, routines  
included in the open-source {\tt Chebfun} package \cite{Cf} were conveniently used 
to determine the Fourier-Legendre coefficients of the functions $f$ and $g$ 
required for the computation of the coefficient matrix $Q_N.$\\

For each of the three types of problems we have studied in the previous section, which give rise to the three algorithms 
we have sketched, we now present the results obtained for two different $g,$ $f$ and  boundary conditions. 
In several tests, we needed an accurate estimate of the exact eigenvalues for the evaluation of the errors in the 
uncorrected and/or corrected numerical ones. In this regard, we decided to consider as ``exact'' 
those provided by a well-established routine or, alternatively, the corrected ones obtained with $N$ very large.  
Further details will be given in the sequel. Before proceeding, we must say that when talking about relative 
errors we usually refer to 
$$\log_{10}\left(|\lambda_{k}^{(N)}-\bar{\lambda}_k|/|\bar{\lambda}_k|\right)  \qquad 
\mbox {or} \qquad \log_{10}\left(|\mu_{k}^{(N)}-\bar{\lambda}_k|/|\bar{\lambda}_k|\right)$$ 
with $\bar{\lambda}_k$ the $k$th reference eigenvalue used.\\

Let us begin with (weakly) regular problems not subject to the Dirichlet boundary condition at $x=-1$ (see Algorithm~\ref{algnod}).
The first results we present confirm that the error in the uncorrected numerical eigenvalues behaves like $O((N+1)^{-p})$
where $p=6-4\ga,$ see Theorem~\ref{tecnod}. In particular, we considered the problems with the following potentials and BCs
\begin{equation}\label{prob1}
q(x) = \cos(2\pi x) +\frac{10 (2-{\rm e}^{-x})}{(1+x)^\ga},  \quad \ga =\frac{1}4,\frac{1}2,\frac{3}4, \quad y(\pm 1) =\pm y'(\pm 1), 
\end{equation}
and we used the classical formula 
\begin{equation}\label{pstima}
p  \approx \log_2\left(\delta\lambda_{k,N}/\delta\lambda_{k,2N+1}\right), \qquad \delta\lambda_{k,N} \equiv |\lambda_k^{(N)}-\lambda_k^{(2N+1)}|,
\end{equation}
for the numerical estimate of the order of convergence. The results we got for the eigenvalues of 
index $k=5,10,20,$ listed in Table~\ref{tab1}, surely confirm the statement of the theorem previously mentioned.

\begin{table}[t]
\caption{Order of convergence for the weakly regular problems (\ref{slp})-(\ref{prob1}).} \label{tab1}
\begin{small}
\begin{center}
\begin{tabular}{|r|c c|c c|c c|}\hline
\multicolumn{7}{|c|}{~}\vspace{-.25cm}\\
\multicolumn{7}{|c|}{$\ga = 0.25, \qquad p = 5$}\\
\multicolumn{7}{|c|}{~}\vspace{-.25cm}\\
\hline
&&&&&&\vspace{-.25cm}\\
$N$ & $\delta\lambda_{5,N}$ & order & $\delta\lambda_{10,N}$ & order & $\delta\lambda_{20,N}$ & order \\
&&&&&&\vspace{-.25cm}\\
\hline
&&&&&&\vspace{-.25cm}\\
 $49$ & $9.9201{\rm E}-08$ & $5.006$ & $1.1937{\rm E}-07$ & $5.002$ & $1.2280{\rm E}-07$ & $4.981$ \\ 
 $99$ & $3.0866{\rm E}-09$ & $5.001$ & $3.7263{\rm E}-09$ & $4.997$ & $3.8890{\rm E}-09$ & $4.988$ \\ 
$199$ & $9.6399{\rm E}-11$ & $5.000$ & $1.1670{\rm E}-10$ & $5.121$ & $1.2255{\rm E}-10$ & $4.865$ \\ 
$399$ & $3.0127{\rm E}-12$ &   --    & $3.3538{\rm E}-12$ &   --    & $4.2064{\rm E}-12$ &   --   \\ 
\hline
\multicolumn{7}{|c|}{~}\vspace{-.25cm}\\
\multicolumn{7}{|c|}{$\ga = 0.5, \qquad p = 4$}\\
\multicolumn{7}{|c|}{~}\vspace{-.25cm}\\
\hline
&&&&&&\vspace{-.25cm}\\
$N$ & $\delta\lambda_{5,N}$ & order & $\delta\lambda_{10,N}$ & order & $\delta\lambda_{20,N}$ & order \\
&&&&&&\vspace{-.25cm}\\
\hline
&&&&&&\vspace{-.25cm}\\
 $49$ & $2.1098{\rm E}-05$ & $4.003$ & $3.0250{\rm E}-05$ & $3.999$ & $3.2895{\rm E}-05$ & $3.981$ \\ 
 $99$ & $1.3159{\rm E}-06$ & $4.001$ & $1.8917{\rm E}-06$ & $4.000$ & $2.0828{\rm E}-06$ & $3.999$ \\ 
$199$ & $8.2192{\rm E}-08$ & $4.000$ & $1.1819{\rm E}-07$ & $4.000$ & $1.3031{\rm E}-07$ & $4.000$ \\ 
$399$ & $5.1362{\rm E}-09$ &   --    & $7.3861{\rm E}-09$ &   --    & $8.1446{\rm E}-09$ &   --   \\ 
\hline
\multicolumn{7}{|c|}{~}\vspace{-.25cm}\\
\multicolumn{7}{|c|}{$\ga = 0.75, \qquad p = 3$}\\
\multicolumn{7}{|c|}{~}\vspace{-.25cm}\\
\hline
&&&&&&\vspace{-.25cm}\\
$N$ & $\delta\lambda_{5,N}$ & order & $\delta\lambda_{10,N}$ & order & $\delta\lambda_{20,N}$ & order \\
&&&&&&\vspace{-.25cm}\\
\hline
&&&&&&\vspace{-.25cm}\\
 $49$ & $1.9714{\rm E}-03$ & $2.999$ & $5.1330{\rm E}-03$ & $2.996$ & $7.5944{\rm E}-03$ & $2.981$ \\ 
 $99$ & $2.4665{\rm E}-04$ & $3.000$ & $6.4360{\rm E}-04$ & $3.000$ & $9.6156{\rm E}-04$ & $2.998$ \\ 
$199$ & $3.0833{\rm E}-05$ & $3.000$ & $8.0475{\rm E}-05$ & $3.000$ & $1.2036{\rm E}-04$ & $3.000$ \\ 
$399$ & $3.8541{\rm E}-06$ &   --    & $1.0060{\rm E}-05$ &   --    & $1.5048{\rm E}-05$ &   --   \\ 
\hline              
\end{tabular}
\end{center}
\end{small}
\end{table}

\begin{figure}[h]
\begin{center}
 \includegraphics[width=12cm, height = 9cm]{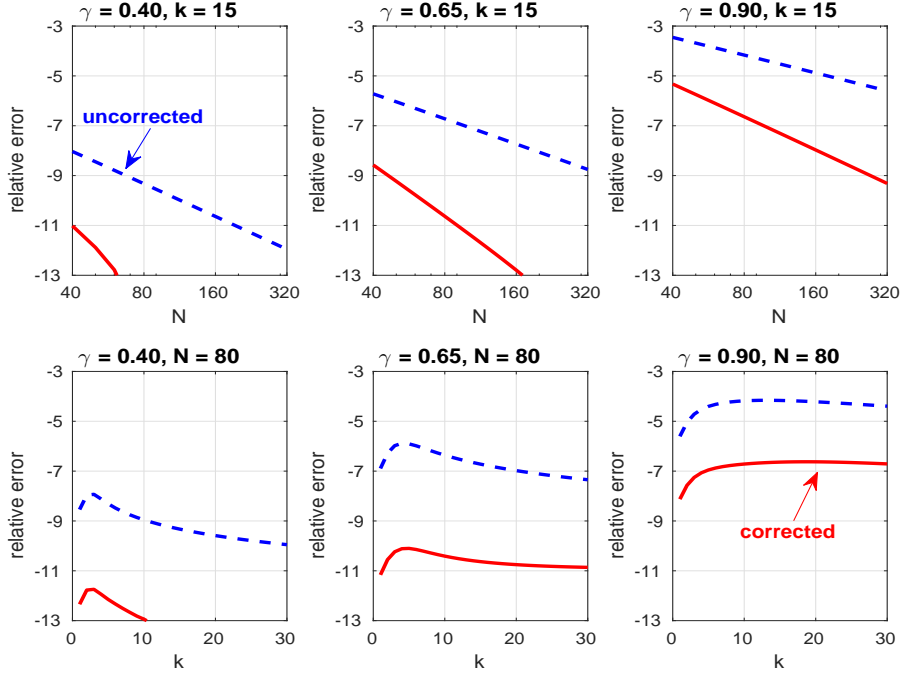}
\caption{Relative errors for the weakly regular problems (\ref{slp})-(\ref{prob2}).}
 \label{fig1}
\end{center}
\end{figure}

Concerning the effectiveness and utility of the application of the a posteriori correction, we applied
Algorithm~\ref{algnod} for solving the problems with  
\begin{equation}\label{prob2}
q(x) = 2x^2 +\frac{5}{((1+x)^2+1)\,(1+x)^\ga}, \quad \ga= 0.4,0.65,0.9, \quad \begin{array}{l} y'(-1) =0, \\ y(1)=0.\end{array}
\end{equation}
For the computation of corresponding reference eigenvalues,
we first of all tried to use the well-established 
and general-purpose codes MATSLISE2 \cite{Led}, SLEDGE \cite{Sledge}, and SLEIGN2 \cite{Sleign} with
a tolerance for the absolute and/or relative error equal to $10^{-13}.$ 
In particular, for the SLEIGN2 routine
we set the input parameter ${\tt NCA}$ equal to two which indicate that the left endpoint is weakly regular.
The approximations we obtained for $\lambda_{15}$ are listed in Table~\ref{tab1_1}. As one can see, 
the number of significant digits for which such estimates agrees decreases as $\ga$ approaches one.
Indeed, this fact is underlined in the documentation of these three softwares and all our tests indicate 
that it is more relevant if $y(-1)\neq 0.$    
We therefore decided 
to use as reference eigenvalues $\lambda_k\approx \bar{\lambda}_k \equiv \mu_k^{(N_T)}$ with $N_T \gg N$ and, in particular, 
for the results shown in Figure~\ref{fig1}, we set $N_T = 3000$ (the  values of $\mu_{15}^{(3000)}$ 
are listed in the  last column of Table~\ref{tab1_1}). In more details, 
in the three subplots at the top of  Figure~\ref{fig1}, the resulting relative errors in the approximation of 
the fifteenth eigenvalue are plotted versus $N$ with $N$ ranging from $40$ to $320.$
For the subplots at the bottom, instead, we fixed $N=80$ and we depict 
the errors for the index $k$ ranging from $1$ to $30.$ The legend of each graphic and of the subsequent ones 
is dashed line and solid line for the errors in the uncorrected numerical eigenvalues and in the corrected ones, respectively.
As one can see,  the correction  improves noticeably the accuracy of the numerical eigenvalues.
As a matter of fact, see the subplots at the top of Figure~\ref{fig1},
it results always 
\begin{equation}\label{guad}
|\mu_{k}^{(N)} - \bar{\lambda}_{k}|\ll |\lambda_{k}^{(2N)} - \bar{\lambda}_{k}|,
\end{equation}
with $k=15.$ On the other hand, from the subplots at the bottom  one deduces that for $N=80$ 
and $1\leq k \leq 30$ the gain resulting from the correction  is at least 
of two significant digits for each eigenvalue but it is frequently much larger.\\

\begin{table}[t]
\caption{Numerical approximations of $\lambda_{15}$ for problems (\ref{slp})-(\ref{prob2}).}\label{tab1_1}
\begin{small}
\begin{center}
\begin{tabular}[h]{|c|ccc|c|}
\hline
\vspace{-.2cm}&&&&\\
$\ga$ & MATSLISE2 \cite{Led} & SLEDGE \cite{Sledge} & SLEIGN2 \cite{Sleign} & $\mu_{15}^{(3000)}$\\
&&&&\vspace{-.25cm}\\
\hline
&&&&\vspace{-.25cm}\\
$0.40$ & $523.9182398826$ & $523.9182763992$ & $523.9182711601$ & $523.9182763990$\\
$0.65$ & $528.1784764701$ & $528.1830115554$ & $528.1791421864$ & $528.1830147149$\\
$0.90$ & $551.0460761225$ & $551.9133020569$ & $550.5488188848$ & $552.2447514722$\\
\hline
\end{tabular}
\end{center}
\end{small}
\end{table}

\begin{figure}[t]
\begin{center}
\includegraphics[width=12cm, height = 9cm]{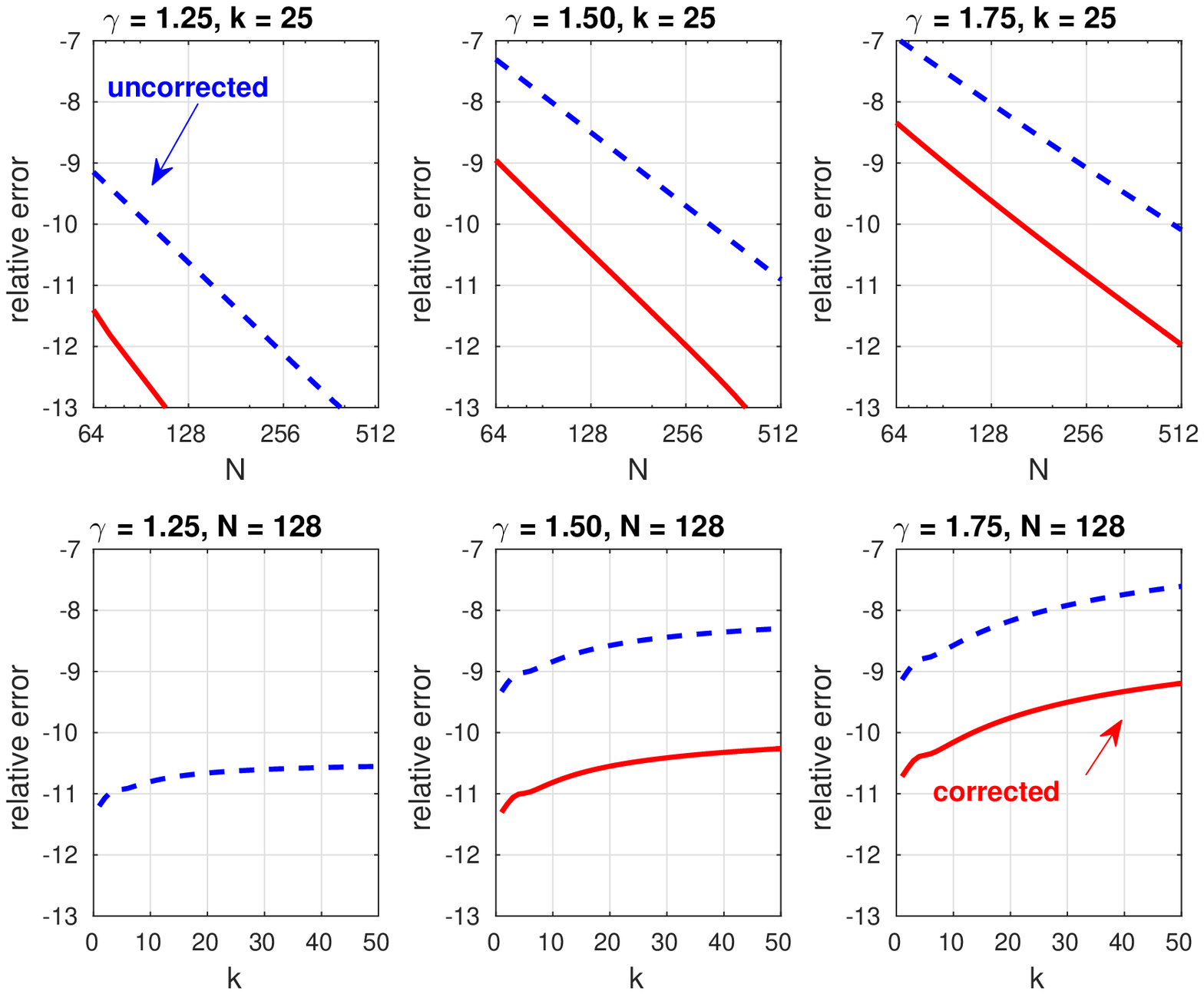}
 \caption{Relative errors for problems (\ref{slp})-(\ref{prob3}).}\label{fig2}
\end{center}
\end{figure}

\begin{figure}[h]
\begin{center}
 \includegraphics[width=12cm, height = 9cm]{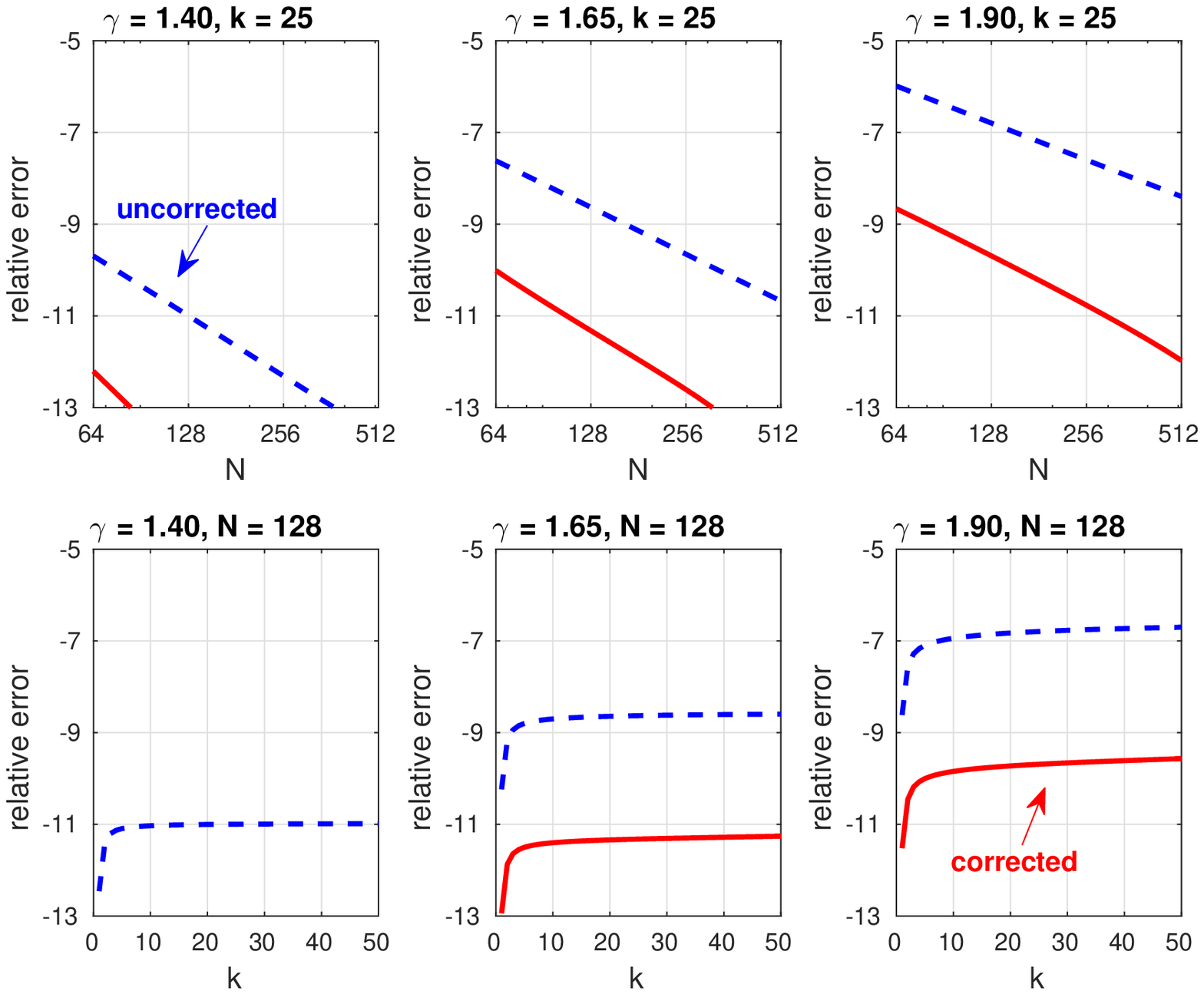}
\caption{Relative errors for problems (\ref{slp})-(\ref{prob4}).}
 \label{fig3}
\end{center}
\end{figure}

The next examples regards problems with $\ga \in (1,2)$ subject to $y(-1)=0$ at the singular endpoint 
which is of type limit-circle. In this case, we applied Algorithm~\ref{algdiri} for getting approximations of the
eigenvalues and we used as reference ``exact'' ones, i.e. $\bar{\lambda}_k,$ those provided by MATSLISE2 \cite{Led}
with a tolerance for the absolute and/or relative errors equal to $10^{-13}.$
It must be said that this choice was motivated only by the fact that
MATSLISE2 is a Matlab code and that the results we are going to present
would have been essentially the same if we had decided to use one of the other two previously mentioned codes.
With the same notation used for the second example, in Figure~\ref{fig2} we represent the errors for the problems with
\begin{equation}
q(x) = \frac{3\left(x\cos(2\pi x)\right)^2}{(1+x)^\ga}, \quad \ga=1.25,1.5,1.75, \qquad y(-1) = y'(1) = 0, \label{prob3}\\
\end{equation}
 while in Figure~\ref{fig3} those corresponding to
\begin{equation}
q(x) =  2\cosh(x) + \frac{2+x}{(1+3x^2)\,(1+x)^\ga}, \,\, \ga=1.4,1.65,1.9, \quad \begin{array}{l}y(-1)=0,\\ y(1)=y'(1).\end{array} \label{prob4}
\end{equation}
For both these examples,  we observe that the spectral matrix 
method already provides accurate approximation for the smallest $\gamma$'s and the application of the correction further improves such estimates.
In particular, see the subplots at the bottom left,  
the relative errors in the first fifty uncorrected numerical eigenvalues, determined with $N=128,$ 
are smaller than $10^{-10}$ while those in the corresponding corrected eigenvalues are smaller than $10^{-13},$
i.e. smaller than the tolerance used for the computation of the reference ones. This is the reason for which 
they are not depicted. Concerning the results 
obtained for $\ga \in [1.5,2),$ the advantage arising from the application of the correction
is undeniable since (\ref{guad}), with $k=25,$ holds almost always.\\


\begin{figure}[t]
\begin{center}
 \includegraphics[width=12.5cm, height = 9cm]{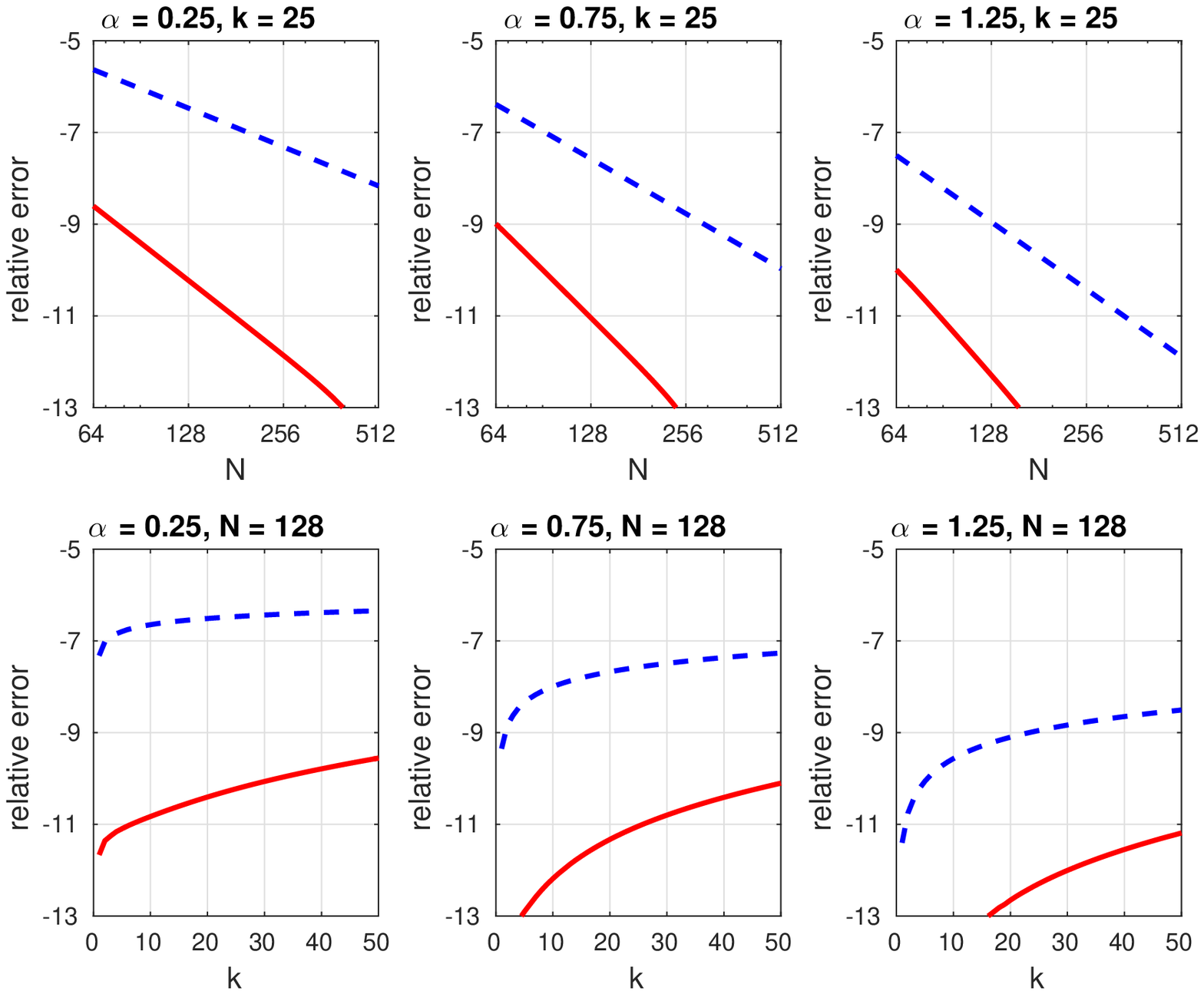}
\caption{Relative errors for problems (\ref{slp})-(\ref{prob6}).} 
 \label{fig4}
\end{center}
\end{figure}

\begin{figure}[t]
\begin{center}
 \includegraphics[width=12cm, height = 9cm]{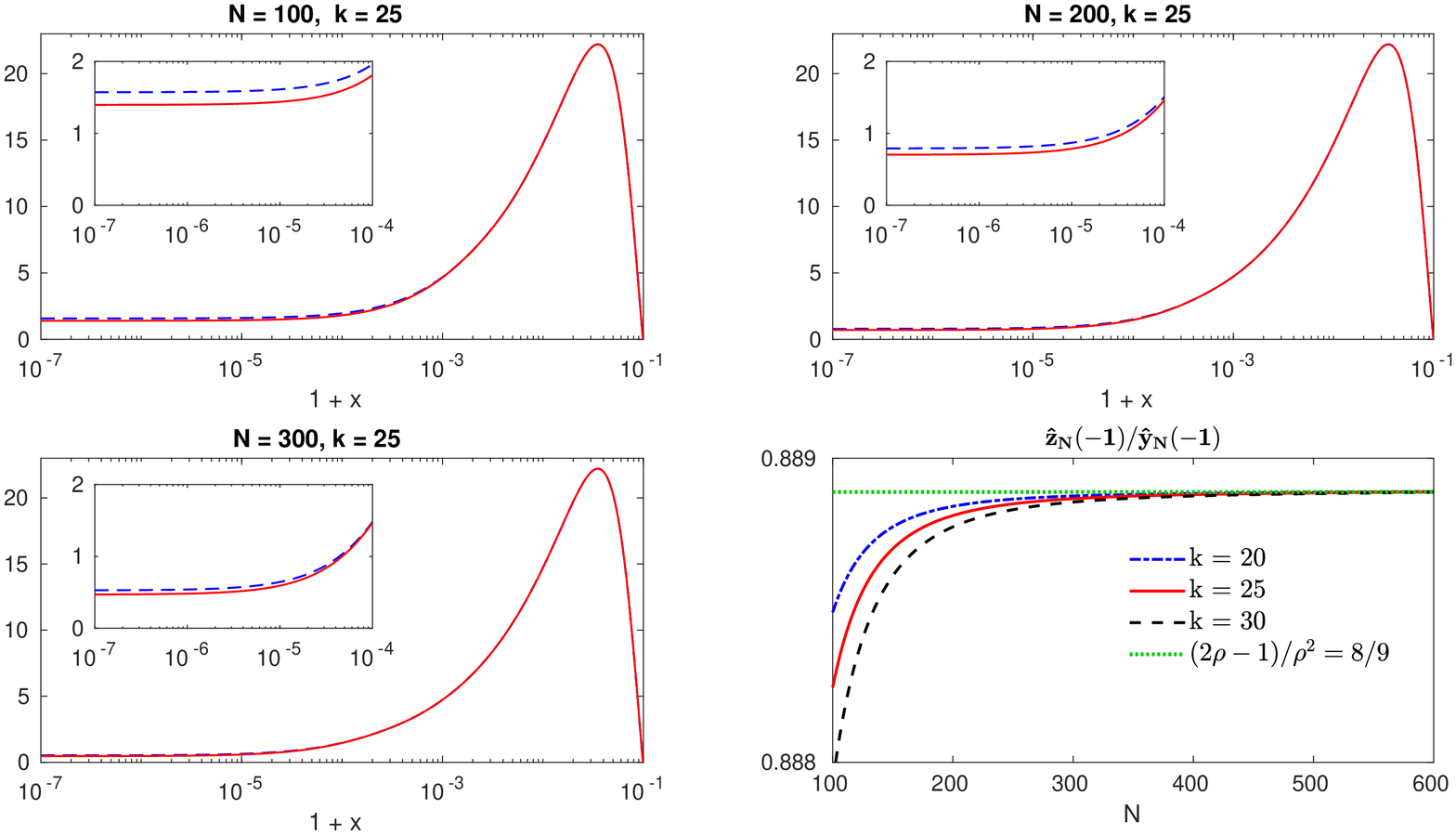}
\caption{Functions $\hat{z}_N(x)$ (solid red line) and $\hat{y}_N(x)$ (dashed blue line) for $x \in (-1,-0.9]$ and 
$\hat{z}_N(-1)/\hat{y}_N(-1)$ versus $N$ for problem (\ref{slp})-(\ref{prob6}) with $\alpha = 0.75.$}  
 \label{fig5}
\end{center}
\end{figure}

Finally, we consider problems with $\ga=2$ and $g(-1)>0.$ In Table~\ref{tab2}, 
we list the experimental orders of convergence, see (\ref{pstima}), 
for the problems with 
\begin{equation}\label{prob5}
q(x) = \log(3+x) + \frac{\alpha\,\cos(4\pi x)}{(1+x)^2}, \quad \alpha =\frac{1}8,\frac{1}2,1,  \qquad y(\pm 1) =0.
\end{equation}
The results obtained are in perfect agreement with the statement of 
Theorem~\ref{teoconvg2}. In fact, $g(x) = \alpha \cos(4\pi x)$ and 
$p=2\sqrt{1+4g(-1)} = 2\sqrt{1+4\alpha},$ see (\ref{ordg2}).\\
Regarding Algorithm~\ref{algg2}, we applied it for solving the 
problems with
\begin{equation}\label{prob6}
q(x) = \frac{1}{1+25x^2} + \frac{\alpha \left(1 + \sinh(1+x)\right)}{(1+x)^2}, \,\,\, \alpha=\frac{1}4,\frac{3}4,\frac{5}4,\quad 
\begin{array}{l} y(-1)=0,\\ y(1) = 2y'(1). \end{array}
\end{equation}
The corresponding errors 
with respect to the eigenvalue estimates provided by MATSLISE2 are represented in Figure~\ref{fig4}.
As one can see, the comments we have done for the previous examples 
concerning the effectiveness of the a posteriori correction surely apply even to this last one.
In addition, in support of (\ref{hzhyn})-(\ref{kapcong}), 
in Figure~\ref{fig5} some graphs of $\hat{z}_N(x)$ and of $\hat{y}_N(x)$ 
(obtained as partial sum  of $\hat{z}_{8000}(x)$) in proximity of $x=-1$ are reported. 
In the same figure, some ratios $\hat{z}_N(-1)/\hat{y}_N(-1)$ are also depicted  
which show that such values approach $(2\rho-1)/\rho^2$ as $N$ increases as stated in 
(\ref{kapcong}) (see also (\ref{ro})).

\begin{table}[ht]
\caption{Order of convergence for problems (\ref{slp})-(\ref{prob5}).} \label{tab2}
\begin{small}
\begin{center}
\begin{tabular}{|r|c c|c c|c c|}
\hline
\multicolumn{7}{|c|}{~}\vspace{-.25cm}\\
\multicolumn{7}{|c|}{$\alpha = 1/8, \qquad p = \sqrt{1+4\alpha} = 2\sqrt{1.5}\approx 2.450$}\\
\multicolumn{7}{|c|}{~}\vspace{-.25cm}\\
\hline
&&&&&&\vspace{-.25cm}\\
$N$ & $\delta\lambda_{5,N}$ & order & $\delta\lambda_{10,N}$ & order & $\delta\lambda_{20,N}$ & order \\
&&&&&&\vspace{-.25cm}\\
\hline
&&&&&&\vspace{-.25cm}\\
$49$  & $1.4443{\rm E}-04$ & $2.448$ & $6.2160{\rm E}-04$ & $2.445$ & $2.8090{\rm E}-03$ & $2.431$ \\ 
$99$  & $2.6461{\rm E}-05$ & $2.449$ & $1.1412{\rm E}-04$ & $2.449$ & $5.2076{\rm E}-04$ & $2.448$ \\ 
$199$ & $4.8448{\rm E}-06$ & $2.449$ & $2.0900{\rm E}-05$ & $2.449$ & $9.5467{\rm E}-05$ & $2.449$ \\ 
$399$ & $8.8697{\rm E}-07$ &   --    & $3.8264{\rm E}-06$ &   --    & $1.7481{\rm E}-05$ &   --    \\ 
\hline
\multicolumn{7}{|c|}{~}\vspace{-.25cm}\\
\multicolumn{7}{|c|}{$\alpha = 1/2, \qquad p =\sqrt{1+4\alpha} =  2\sqrt{3}\approx 3.464$}\\
\multicolumn{7}{|c|}{~}\vspace{-.25cm}\\
\hline
&&&&&&\vspace{-.25cm}\\
$N$ & $\delta\lambda_{5,N}$ & order & $\delta\lambda_{10,N}$ & order & $\delta\lambda_{20,N}$ & order \\
&&&&&&\vspace{-.25cm}\\
\hline
&&&&&&\vspace{-.25cm}\\
$49$  & $8.4050{\rm E}-05$ & $3.463$ & $4.0854{\rm E}-04$ & $3.459$ & $2.4019{\rm E}-03$ & $3.438$ \\ 
$99$  & $7.6244{\rm E}-06$ & $3.464$ & $3.7163{\rm E}-05$ & $3.463$ & $2.2161{\rm E}-04$ & $3.462$ \\ 
$199$ & $6.9096{\rm E}-07$ & $3.464$ & $3.3691{\rm E}-06$ & $3.464$ & $2.0117{\rm E}-05$ & $3.464$ \\ 
$399$ & $6.2613{\rm E}-08$ &   --    & $3.0530{\rm E}-07$ &   --    & $1.8233{\rm E}-06$ &   --    \\ 
\hline
\multicolumn{7}{|c|}{~}\vspace{-.25cm}\\
\multicolumn{7}{|c|}{$\alpha = 1, \qquad p = 2\sqrt{1+4\alpha}=2\sqrt{5}\approx 4.472$}\\
\multicolumn{7}{|c|}{~}\vspace{-.25cm}\\
\hline
&&&&&&\vspace{-.25cm}\\
$N$ & $\delta\lambda_{5,N}$ & order & $\delta\lambda_{10,N}$ & order & $\delta\lambda_{20,N}$ & order \\
&&&&&&\vspace{-.25cm}\\
\hline
&&&&&&\vspace{-.25cm}\\
$49$  & $8.6382{\rm E}-06$ & $4.470$ & $4.5493{\rm E}-05$ & $4.465$ & $3.2872{\rm E}-04$ & $4.425$ \\ 
$99$  & $3.8980{\rm E}-07$ & $4.472$ & $2.0601{\rm E}-06$ & $4.471$ & $1.5299{\rm E}-05$ & $4.469$ \\ 
$199$ & $1.7566{\rm E}-08$ & $4.472$ & $9.2871{\rm E}-08$ & $4.472$ & $6.9082{\rm E}-07$ & $4.472$ \\ 
$399$ & $7.9135{\rm E}-10$ &   --    & $4.1846{\rm E}-09$ &   --    & $3.1136{\rm E}-08$ &   --   \\ 
\hline      
\end{tabular}
\end{center}
\end{small}
\end{table}

\section*{Acknowledgments}
The author is very indebted to Prof.\,Paolo Ghelardoni for helpful discussions and suggestions.

\end{document}